\documentclass[11pt,a4paper,reqno]{amsart}
\usepackage{amsmath,amssymb,amsthm,mathrsfs,hyperref}
\usepackage[margin=30mm]{geometry}
\hypersetup{colorlinks=true,citecolor=blue,linkcolor=blue,urlcolor=blue}
\numberwithin{equation}{section}
\newtheorem{theorem}{Theorem}[section]
\newtheorem{proposition}[theorem]{Proposition}
\newtheorem{lemma}[theorem]{Lemma}
\newtheorem{corollary}[theorem]{Corollary}
\newtheorem{definition}[theorem]{Definition}
\theoremstyle{remark}\newtheorem{remark}[theorem]{Remark}
\newcommand{\R}{\mathbb R}\newcommand{\C}{\mathbb C}\newcommand{\Z}{\mathbb Z}\newcommand{\N}{\mathbb N}\newcommand{\Sspace}{\mathcal S}\newcommand{\supp}{\operatorname{supp}}
\title[Critical Gaussian asymptotics and STFT consequences]{Critical Gaussian Asymptotics for Directional Superoscillations: STFT and Gabor Consequences}
\author{Younis Ahmad Bhat}\author{Sunit Ghosh}\author{Jitendriya Swain}
\address{Department of Mathematics, Indian Institute of Technology Guwahati, Guwahati 781039, India}
\email{younis@iitg.ac.in}\email{g.sunit@iitg.ac.in}\email{jitumath@iitg.ac.in}
\subjclass[2020]{42B10, 42C15, 46E10, 35Q41}
\keywords{Superoscillation, short-time Fourier transform, entire functions, Gabor frames, higher-order Schr\"odinger equation, uniform electric field}
\begin{document}\begin{abstract}
We study directional Aharonov--Berry superoscillations in $\R^d$ with widening anisotropic Gaussian envelopes. We prove that the rescaled total $L^2$ mass has a sharp transition governed by the matrix $D=Q-(a^2-1)vv^T$. At criticality, $D$ has a one-dimensional kernel, and the mass obeys an exact $N^{d/2+1/4}$ law with an explicit Gamma-function constant. After anisotropic rescaling, the intensity converges globally in $L^1$ to a transverse Gaussian times a longitudinal quartic exponential. Moyal's identity gives the corresponding total STFT-energy law, but does not imply a local phase-space concentration theorem. We compare these global statements with the earlier quartic critical-window analysis of Berry and Shukla for periodic superoscillations. Supporting results include quantitative fixed-Gaussian STFT expansions, entire-growth convergence with a critical-step obstruction, and conditional Gabor-frame and dual-window expansions. Weighted STFT estimates for polynomial electric-field evolution are obtained as applications. No claim of priority is made for quartic criticality, qualitative STFT supershift preservation, or polynomial electric-field persistence as such.
\end{abstract}\maketitle\enlargethispage{3pt}
\section{Introduction}
Superoscillations arise when finite combinations of bounded frequencies approximate a wave with a frequency outside their spectral range. The Aharonov--Berry sequence
\begin{equation}\label{intro:FN}
 F_N(\theta,a)=\left(\cos\frac{\theta}{N}+ia\sin\frac{\theta}{N}\right)^N
 =\sum_{j=0}^N c_{j,N}(a)e^{is_{j,N}\theta},\quad
 s_{j,N}=1-\frac{2j}{N},\quad a>1,
\end{equation}
where $c_{j,N}(a)=\binom Nj ((1+a)/2)^{N-j}((1-a)/2)^j$, converges to $e^{ia\theta}$ locally uniformly. Quantitative information becomes substantially more delicate when the observation variable grows with $N$ or when the signals are transformed or propagated.

Critical windowing of superoscillations has an earlier precedent. Berry and Shukla \cite{BerryShukla} analyze a periodic canonical superoscillatory function with a periodic Gaussian-type window. Their Section~6 identifies a quartic critical maximum at the window transition (equation~(6.8)) and Gamma-function expressions in the associated uncertainty-product asymptotics (equation~(6.9)). We therefore make no priority claim for the occurrence of quartic criticality or Gamma functions in windowed superoscillations. Here the different object is the full $L^2(\R^d)$ mass for the nonperiodic directional family with arbitrary $Q>0$, including the matrix criterion, the exact critical mass constant, and global anisotropic $L^1$ convergence of intensity. The latter is a global dominated-convergence theorem, not solely a local quartic approximation.

The short-time Fourier transform (STFT) makes this problem a question about phase-space concentration. Alpay, De Martino, Diki and Struppa \cite{ADDS} investigated STFTs and Zak transforms of superoscillations, including Gaussian and Hermite windows, Gabor frames and evolution. Alpay, Colombo, Diki, Sabadini and Struppa \cite{ACDSS} compute the Segal--Bargmann transforms of Gaussian-localized superoscillatory sequences as coherent-state superpositions; consequently, Gaussian-transform identities themselves are not a priority claim here. Analytic supershifts and persistence for Schr\"odinger-type equations, including higher-order operators and uniform electric fields, have also been studied extensively; see \cite{AMS,ACSST,Electric}. Our purpose is not to reprove qualitative persistence. Instead we obtain explicit norm expansions and precisely formulated frame perturbation statements, and isolate a critical Gaussian scaling that demands a different limiting analysis.

A closely related multivariable supershift theory was developed by Colombo, Pinton, Sabadini and Struppa in 2023, and subsequent analyticity refinements appeared in 2025; these works already treat analytic superpositions and infinite-order differential operators. Hence quantitative remainder estimates in a new norm, not multivariable analyticity alone, must provide the novelty of the present work \cite{CPSS,CSSY}. The theorem at the heart of this article describes the boundary
$D=Q-(a^2-1)vv^T\ge0$ of the Gaussian localization regime. At
$(a^2-1)v^TQ^{-1}v=1$, exactly one quadratic direction degenerates, and
we prove that the squared mass is asymptotic to an explicit positive
constant times $N^{d/2+1/4}$. Moreover the anisotropically rescaled
intensities converge globally in $L^1$ to a transverse Gaussian times a
longitudinal quartic exponential. This exact constant and global limit,
not merely the known qualitative STFT supershift phenomenon, are the
principal asymptotic conclusions. The corresponding STFT result is a
total-energy identity obtained by Moyal's theorem, not a pointwise
phase-space localization statement.

Three technical distinctions guide the paper. First, an unwindowed plane wave is not in $L^2(\R^d)$, whereas multiplying by a fixed Gaussian makes Hilbert-space and global STFT estimates possible. Second, nonlinear frequency maps are compatible with compactly supported smooth windows; such localization is important because the $\ell^1$-norm of the signed coefficients in \eqref{intro:FN} equals $a^N$. Third, the exact holomorphic growth of a Gaussian STFT does not, on its own, imply convergence in a fixed entire-function Banach step. We distinguish pointwise kernel membership, convergence in larger fixed growth steps, and failure at the endpoint. For context, convergence of the classical sequence in order-one entire-function growth spaces has previously been studied in \cite{ACS2021}; the Gaussian STFT growth topology used here is different.

\subsection*{Positioning relative to previous work}
Alpay--De Martino--Diki--Struppa \cite{ADDS} established qualitative STFT supershift preservation and considered Zak transforms, Gabor frames, Gaussian and Hermite windows, and an evolution problem. These topics themselves are not new. Our fixed-window statements provide explicit second-order remainders in global norms; the frame-operator and canonical-dual expansions quantify reconstruction under a separately assumed frame property. In contrast to a blanket inference from a nonzero window to a uniform Zak lower bound, Proposition~\ref{prop:zak} records a direct failure example. Our global widening-Gaussian dichotomy, critical mass rate, exact constant, and anisotropic $L^1$ profile are the central threshold results of this paper; the corresponding STFT statement concerns total energy, not local phase-space concentration.

The multivariable analytic supershift framework of \cite{CPSS,CSSY} and the earlier analytic-convolution theory \cite{AMS,Aoki} already handle broad classes of parameter-dependent entire functions. We therefore do not claim novelty for general multivariable supershifts or infinite-order operators. Likewise, even-polynomial momentum Hamiltonians in uniform electric fields appear in \cite{Electric}; our evolution section provides a quantitatively weighted STFT application rather than a new qualitative persistence theorem. The simple $L^2$-based persistence estimate is a consequence of unitarity, while the finite-time polynomially weighted version requires additional estimates on Schwartz seminorms.

\subsection*{Precise scope of the critical-threshold contribution}
The periodic model of \cite{BerryShukla} uses a canonical superoscillatory
function on the circle, multiplied by a periodic Gaussian-type window whose
width is varied. At its critical window width, their equation~(6.8)
contains a quartic local maximum and equation~(6.9) evaluates uncertainty
quantities involving Gamma-function constants. Thus the cancellation of a
quadratic term and the presence of a quartic profile are known mechanisms.
Our family differs in the scaling and the global observable:
\[
 h_N(x)=F_N(v\cdot x,a)\exp\!\left(-\frac{x^TQx}{2N}\right),
 \qquad x\in\mathbb R^d,\quad Q=Q^T>0.
\]
Theorem~\ref{thm:globalcritical} characterizes finite versus divergent
rescaled total mass by the rank-one matrix $Q-(a^2-1)vv^T$.
Theorem~\ref{thm:quarticcritical} then gives, at the boundary,
(i) the exact mass exponent and leading constant and (ii) strong global
$L^1$ convergence of anisotropically rescaled intensities, including
transverse Gaussian directions. These are the precise claims distinguished from the cited precedents. The reviewed works \cite{BerryShukla,ACDSS,ADDS} do not state
this combined global anisotropic theorem in the respective settings
examined here; this is a comparison of their stated results, not a claim of
exhaustive priority across all asymptotic literature.

\subsection*{Analytic difficulties and organization}
Three issues motivate our proof structure. The signed Bernstein coefficients have $\ell^1$ norm $a^N$, so compact convergence alone cannot be used for global STFT norms. At critical Gaussian scaling, the quadratic term cancels along one direction; a quartic Taylor coefficient yields the correct rate but global domination is indispensable for the leading constant. Finally, frame perturbation statements require an actual lower frame bound for the limiting Gabor system, and the ordinary STFT need not follow a single classical trajectory for nonlinear polynomial dispersion.

We first establish estimates for fixed Gaussian envelopes. Next come
the global Gaussian threshold and its sharp quartic boundary law.
We then investigate entire growth and Gabor reconstruction. The last
sections concern electric-field evolution and co-moving observables.
Quantitative coupled-frequency supershifts and localization along a
single nonlinear trajectory are not claimed here.

\section{Notation and elementary bounds}
We use the angular-frequency Fourier transform
\[\widehat f(\xi)=\int_{\R^d}f(x)e^{-i\xi\cdot x}\,dx,\qquad f(x)=(2\pi)^{-d}\int_{\R^d}\widehat f(\xi)e^{i\xi\cdot x}\,d\xi.\]
For $h\in L^2(\R^d)$, its STFT is
\begin{equation}\label{STFT} V_hf(u,\omega)=\int_{\R^d}f(x)\overline{h(x-u)}e^{-i\omega\cdot x}\,dx.\end{equation}
The time and frequency shifts are $T_uf(x)=f(x-u)$ and $M_\omega f(x)=e^{i\omega\cdot x}f(x)$. Moyal's identity and Cauchy--Schwarz imply
\begin{equation}\label{STFTbounds}
 \|V_hf\|_2=(2\pi)^{d/2}\|h\|_2\|f\|_2,\qquad
 \|V_hf\|_\infty\le\|h\|_2\|f\|_2.
\end{equation}
Interpolation yields, for $2\le q\le\infty$,
\begin{equation}\label{qbound}
 \|V_hf\|_q\le (2\pi)^{d/q}\|h\|_2\|f\|_2.
\end{equation}
We use the $M^1$ modulation space associated with a fixed nonzero Schwartz window $\phi$:
\[\|f\|_{M^1}=\int_{\R^{2d}}|V_\phi f(u,\omega)|\,du\,d\omega.\]
Different Schwartz windows give equivalent norms. For a full-rank time--frequency lattice $\Gamma\subset\R^{2d}$, write
$\pi(u,\omega)=M_\omega T_u$ and set
\[
 C_wf=(\langle f,\pi(\gamma)w\rangle)_{\gamma\in\Gamma},
 \qquad D_wc=\sum_{\gamma\in\Gamma}c_\gamma\pi(\gamma)w,
 \qquad S_w=D_wC_w,
\]
where $\langle f,g\rangle=\int f\overline g$ is linear in its
first entry. The synthesis sum is interpreted as an $L^2$-limit whenever
$w\in M^1$, and the same convention is used for an arbitrary frame window. We recall the standard sampling consequence of $M^1$ theory \cite{Gro}: for a fixed lattice there is $K_\Gamma>0$ such that
\begin{equation}\label{M1frame}
\|C_w\|_{L^2\to\ell^2(\Gamma)}=\|D_w\|_{\ell^2(\Gamma)\to L^2}\le K_\Gamma\|w\|_{M^1}.
\qquad\text{Hence }\|C_w\|+\|D_w\|\le2K_\Gamma\|w\|_{M^1}.
\end{equation}
For $B>0$ and $p>0$, let
\[A_{p,B}(\C^d)=\left\{H\in\mathcal O(\C^d):\sup_{z\in\C^d}|H(z)|e^{-B|z|^p}<\infty\right\}.\]
We write explicitly
\[
 \|H\|_{p,B}=\sup_{z\in\C^d}|H(z)|e^{-B|z|^p}.
\]
We do not identify this entire-growth scale with the distinct Hermite
spectral coefficient spaces sometimes also denoted $A_p$.

\section{Global Gaussian-localized supershifts}
Fix $v\in\R^d\setminus\{0\}$, a symmetric positive-definite $d\times d$ matrix $Q$, and a centre $x_0\in\R^d$. Write
\[g(x)=e^{-\frac12(x-x_0)^TQ(x-x_0)},\quad
 f_N(x)=F_N(v\cdot x,a)g(x),\quad f_a(x)=e^{ia v\cdot x}g(x).\]
The choice of a single direction $v$ permits the global estimates below. Quantitative estimates for genuinely coupled nonlinear frequency curves are not claimed; see the final discussion.

\begin{lemma}[Global expansion]\label{lem:global}
For each integer $m\ge 1$ there are polynomials $P_0,\ldots,P_m$ in one real variable, with $P_0=1$ and $P_1(\theta)=(a^2-1)\theta^2/2$, such that
\begin{equation}\label{globalexp}
\left\|f_N-e^{ia v\cdot x}g(x)\sum_{j=0}^mN^{-j}P_j(v\cdot x)\right\|_{L^2(\R^d)}=O(N^{-m-1}).
\end{equation}
In particular,
\[P_2(\theta)=-\frac{ia(a^2-1)}3\theta^3+\frac{(a^2-1)^2}{8}\theta^4.\]
\end{lemma}
\begin{proof}
Put $\theta=v\cdot x$ and $L(z)=\log(\cos z+ia\sin z)$, where the
branch is chosen near the origin with $L(0)=0$. Its first derivatives are
\[
 L'(0)=ia,\qquad L''(0)=a^2-1,\qquad
 L^{(3)}(0)=-2ia(a^2-1).
\]
Choose $\delta>0$ such that $L$ is analytic on $|z|\le2\delta$.
Taylor's theorem with Cauchy's estimate supplies, for each $m$, a
constant $C_m$ such that
\[
 L(z)=iaz+\sum_{k=2}^{m+2}\ell_kz^k+R_{m+2}(z),
 \qquad |R_{m+2}(z)|\le C_m|z|^{m+3}
 \quad (|z|\le\delta).
\]
Hence
\[
 NL(\theta/N)-ia\theta
 =\sum_{j=1}^{m+1}N^{-j}\ell_{j+1}\theta^{j+1}
 +O_m(N^{-m-2}|\theta|^{m+3}).
\]
The coefficients $P_j(\theta)$ are obtained by expanding the exponential
of the right-hand side in powers of $N^{-1}$: explicitly
\[
 P_j(\theta)=\sum_{r=1}^{j}\frac1{r!}
 \sum_{\substack{k_1,\ldots,k_r\ge1\\ k_1+\cdots+k_r=j}}
 \prod_{\nu=1}^r\ell_{k_\nu+1}\theta^{k_\nu+1}
 \quad(j\ge1),\qquad P_0=1.
\]
This also verifies the stated values of $P_1$ and $P_2$.

For completeness we justify the uniform remainder used under the integral.
Fix $m$, and choose $\varepsilon>0$ sufficiently small that
$\varepsilon<1/4$ and $(m+3)\varepsilon<1$. In the central region
$|\theta|\le N^\varepsilon$, the exponent $NL(\theta/N)-ia\theta$
is bounded, uniformly for large $N$. Taylor's integral formula for
$e^w$, applied to the displayed finite sum, and the residual bound
above imply
\[
 \left|F_N(\theta,a)-e^{ia\theta}
 \sum_{j=0}^mN^{-j}P_j(\theta)\right|
 \le C_mN^{-m-1}(1+|\theta|^{K_m})
 \qquad(|\theta|\le N^\varepsilon),
\]
for an integer $K_m$ depending only on $m$. Here is a uniform justification that does not treat the exponential as
a finite polynomial. Set
\[
 W_N(\theta)=\sum_{j=1}^{m+1}N^{-j}\ell_{j+1}\theta^{j+1}.
\]
For $|\theta|\le N^\varepsilon$ and sufficiently large $N$,
$|\theta|/N\le1$ and
$|W_N(\theta)|\le C_m\theta^2/N$. The exponential Taylor
formula with integral remainder gives
\[
 \left|e^{W_N(\theta)}-
 \sum_{r=0}^{m}\frac{W_N(\theta)^r}{r!}\right|
 \le C_m|W_N(\theta)|^{m+1}e^{|W_N(\theta)|}
 \le C_mN^{-m-1}|\theta|^{2m+2}e^{C_m\theta^2/N}.
\]
The displayed sum over $r$ is genuinely finite. Expanding its
polynomial factors in $N^{-1}$ and collecting all terms of total
$N^{-1}$ degree at most $m$ produces
$\sum_{j=0}^{m}N^{-j}P_j(\theta)$; every uncollected monomial has
degree at least $m+1$, and its finitely many polynomial factors are
bounded by $N^{-m-1}(1+|\theta|^{K_m})$ for $N\ge1$.
The logarithmic Taylor remainder has modulus at most
$C_mN^{-m-2}|\theta|^{m+3}$. The inequality
$|e^{z+w}-e^z|\le |w|e^{|z|+|w|}$ shows that its
contribution is absorbed into the same polynomial remainder.
Finally,
\[
 \sup_{|\theta|\le N^\varepsilon}e^{C_m\theta^2/N}
 \le e^{C_mN^{2\varepsilon-1}}\le e^{C_m},
\]
so all constants are independent of $N$.

For every real $\theta$,
\[
 |F_N(\theta,a)|^2
 =\bigl(1+(a^2-1)\sin^2(\theta/N)\bigr)^N
 \le \exp\bigl((a^2-1)\theta^2/N\bigr).
\]
Since $Q$ is positive definite, there are constants $c_0,C_0>0$ such that
$|g(x)|\le C_0e^{-c_0|x|^2}$. By the exact amplitude estimate and
$|\theta|\le|v||x|$, choose $N_0$ so large that
$(a^2-1)|v|^2/(2N)\le c_0/2$ for $N\ge N_0$. Then
$|F_N(\theta,a)g(x)|\le C_0e^{-c_0|x|^2/2}$.
The same Gaussian majorant with an additional fixed polynomial applies
to each term of the truncated expansion. On
$\{|v\cdot x|>N^\varepsilon\}$, the inequality
$|x|\ge N^\varepsilon/|v|$ and Gaussian tail integration give
$O(e^{-c'N^{2\varepsilon}})$ in $L^2$ (after changing $c'$).
On the complementary region the preceding polynomial remainder,
multiplied by $g$, has $L^2$ norm at most
\[
 C_mN^{-m-1}\bigl\|(1+|v\cdot x|^{K_m})g\bigr\|_2.
\]
Adding the two regions proves the assertion.\end{proof}

\begin{theorem}[Quantitative multidimensional STFT supershift]\label{thm:STFT}
For $h\in L^2(\R^d)\setminus\{0\}$ and $2\le q\le\infty$, put $r_1(x)=\frac{a^2-1}{2}(v\cdot x)^2f_a(x)$ and
\[r_2(x)=\left[-\frac{ia(a^2-1)}3(v\cdot x)^3+\frac{(a^2-1)^2}{8}(v\cdot x)^4\right]f_a(x).\]
Then
\begin{equation}\label{stftexpansion}
\|V_hf_N-V_hf_a-N^{-1}V_hr_1-N^{-2}V_hr_2\|_{L^q(\R^{2d})}=O(N^{-3}).
\end{equation}
In particular $N(V_hf_N-V_hf_a)\to V_hr_1$ in $L^q$, and
\[\lim_{N\to\infty}N\|V_hf_N-V_hf_a\|_2=(2\pi)^{d/2}\|h\|_2\|r_1\|_2>0.\]
\end{theorem}
\begin{proof}
Lemma~\ref{lem:global} with $m=2$ shows that
$E_N=f_N-f_a-N^{-1}r_1-N^{-2}r_2$ satisfies
$\|E_N\|_2\le CN^{-3}$. Since the STFT is linear,
\[
 V_hE_N=V_hf_N-V_hf_a-N^{-1}V_hr_1-N^{-2}V_hr_2.
\]
The interpolation estimate~\eqref{qbound} gives
$\|V_hE_N\|_q\le(2\pi)^{d/q}\|h\|_2CN^{-3}$.
Multiplying the identity by $N$ and taking the limit proves strong
$L^q$ convergence because $\|N^{-1}V_hr_2+NV_hE_N\|_q\to0$.
For $q=2$, Moyal's equality gives
$\|V_hr_1\|_2=(2\pi)^{d/2}\|h\|_2\|r_1\|_2$.
Since $v\ne0$ and the Gaussian is positive, $r_1$ is nonzero in
$L^2$, establishing the positive sharp leading coefficient.\end{proof}

\section{A global dichotomy at the critical Gaussian scale}
The pointwise profile of Proposition~\ref{prop:critical} has a stronger global consequence. The statement below detects precisely when the widening Gaussian suppresses the superoscillatory amplification on the parabolic spatial scale.

\begin{theorem}[Sharp global $L^2$ threshold]\label{thm:globalcritical}
Let $a>1$, $v\in\R^d\setminus\{0\}$, and $Q=Q^T>0$. Define
\[
 h_N(x)=F_N(v\cdot x,a)\exp\left(-\frac{x^TQx}{2N}\right),\qquad
 D=Q-(a^2-1)vv^T,
\]
and let
\[
 H_N(y)=e^{-ia\sqrt N\,v\cdot y}h_N(\sqrt N\,y).
\]
Then the following alternatives hold:
\begin{enumerate}
\item If $D>0$, then
\begin{equation}\label{eq:strongcrit}
 H_N\longrightarrow H_\infty(y):=\exp\left(-\frac12y^TDy\right)
 \quad\text{in }L^2(\R^d),
\end{equation}
and
\begin{equation}\label{eq:masscrit}
 N^{-d/2}\|h_N\|_2^2=\|H_N\|_2^2\longrightarrow
 \int_{\R^d}e^{-y^TDy}\,dy=\frac{\pi^{d/2}}{\sqrt{\det D}}.
\end{equation}
\item If $D$ is not positive definite, then
\begin{equation}\label{eq:divergencecrit}
 N^{-d/2}\|h_N\|_2^2=\|H_N\|_2^2\longrightarrow+\infty.
\end{equation}
\end{enumerate}
In particular the transition occurs precisely at
\[
(a^2-1)v^TQ^{-1}v=1.
\]
\end{theorem}
\begin{proof}
Write $c=a^2-1>0$ and $\theta=v\cdot y$. The exact modulus identity
\[
 |F_N(\sqrt N\,\theta,a)|^2
 =\left(1+c\sin^2\frac\theta{\sqrt N}\right)^N
\]
and the inequalities $\log(1+t)\le t$ and $|\sin s|\le|s|$ give
\begin{equation}\label{eq:criticalmajorant}
 |H_N(y)|^2\le e^{-y^TQy+c(v\cdot y)^2}=e^{-y^TDy}
 \quad\text{for all }N\ge1,\ y\in\R^d.
\end{equation}
To identify the phase as well as the modulus, use the analytic logarithm of
$\cos z+ia\sin z$ near the origin. Its Taylor expansion is
\[
 \log(\cos z+ia\sin z)=iaz+\frac c2z^2+O_a(z^3).
\]
For each fixed $y$ this implies
\[
 N\log\left(\cos\frac{v\cdot y}{\sqrt N}
             +ia\sin\frac{v\cdot y}{\sqrt N}\right)
 =ia\sqrt N\,v\cdot y+\frac c2(v\cdot y)^2+O_{a,y}(N^{-1/2}).
\]
Consequently $H_N(y)\to e^{-y^TDy/2}$ pointwise on $\R^d$.

Suppose first that $D>0$. By \eqref{eq:criticalmajorant},
$|H_N(y)|\le e^{-y^TDy/2}$, and also $|H_\infty(y)|=e^{-y^TDy/2}$.
Thus
$|H_N-H_\infty|^2\le4e^{-y^TDy}$,
which is integrable. Dominated convergence proves \eqref{eq:strongcrit};
applying dominated convergence to $|H_N|^2$ proves \eqref{eq:masscrit}.
The Gaussian integral follows by orthogonally diagonalizing $D$.

Now suppose that $D$ is not positive definite. Since $D$ is symmetric,
there exists a nonzero vector $w$ with $w^TDw\le0$.
If $D$ has a negative eigenvalue, then $e^{-y^TDy}$ grows exponentially
along a cone around a negative-eigenvalue direction. If $D$ is
positive semidefinite but singular, an orthogonal coordinate change
shows that $e^{-y^TDy}$ is constant along a nontrivial null direction.
In either case,
\[
 \int_{\R^d}e^{-y^TDy}\,dy=+\infty.
\]
The pointwise convergence $|H_N(y)|^2\to e^{-y^TDy}$ and Fatou's lemma yield
\[
 +\infty=\int_{\R^d}\liminf_{N\to\infty}|H_N(y)|^2\,dy
 \le\liminf_{N\to\infty}\int_{\R^d}|H_N(y)|^2\,dy.
\]
This proves \eqref{eq:divergencecrit}. Finally the rank-one criterion
$D>0\iff c\,v^TQ^{-1}v<1$ follows by factoring
$D=Q^{1/2}(I-cQ^{-1/2}vv^TQ^{-1/2})Q^{1/2}$.
\end{proof}

\begin{corollary}[Critical STFT energy dichotomy]\label{cor:stftcritical}
For every nonzero $k\in L^2(\R^d)$,
\begin{equation}\label{eq:stftcritenergy}
 N^{-d/2}\|V_kh_N\|_{L^2(\R^{2d})}^2
 =(2\pi)^d\|k\|_2^2\|H_N\|_2^2.
\end{equation}
Consequently, if $D>0$ the left-hand side converges to
$(2\pi)^d\|k\|_2^2\pi^{d/2}/\sqrt{\det D}$, and otherwise it diverges
to $+\infty$.
\end{corollary}
\begin{proof}
Moyal's identity gives $\|V_kh_N\|_2^2=(2\pi)^d\|k\|_2^2\|h_N\|_2^2$.
The substitution $x=\sqrt N y$ yields $\|h_N\|_2^2=N^{d/2}\|H_N\|_2^2$ because the demodulation factor has modulus one. Apply Theorem~\ref{thm:globalcritical}.
\end{proof}

\begin{remark}[Meaning of the transition]
The energy dichotomy is global and sharp, but it does not imply a pointwise STFT limit in moving phase-space coordinates. In particular, \eqref{eq:stftcritenergy} is an identity for the total $L^2$ mass, not a concentration theorem for individual time-frequency cells.
\end{remark}

\section{Sharp quartic asymptotics at the critical Gaussian threshold}\label{sec:quarticcritical}
The global dichotomy of Theorem~\ref{thm:globalcritical} leaves open the precise rate of divergence at the boundary. We now determine this rate, its exact leading constant, and the limiting anisotropically rescaled intensity profile. Our argument controls the full integral, not only the Taylor expansion near the origin.

Put $c=a^2-1>0$ and
\[
 \kappa=\frac c3+\frac{c^2}{2}>0,
 \qquad \delta(t)=ct^2-\log(1+c\sin^2t),\quad t\in\R.
\]
The inequality $\sin^2t\le t^2$ implies $\delta(t)\ge0$, with strict inequality when $t\ne0$.

\begin{lemma}[Global coercivity of the superoscillatory deficit]\label{lem:coercivity}
There is a constant $k_c>0$, depending only on $c$, for which
\begin{equation}\label{eq:deficitcoercivity}
 \delta(t)\ge k_c\frac{t^4}{1+t^2},\qquad t\in\R.
\end{equation}
Moreover, $\delta(t)=\kappa t^4+O_c(t^6)$ as $t\to0$.
\end{lemma}
\begin{proof}
Using $\sin^2t=t^2-t^4/3+O(t^6)$ and $\log(1+s)=s-s^2/2+O(s^3)$, we have
\[
\log(1+c\sin^2t)=ct^2-\left(\frac c3+\frac{c^2}{2}\right)t^4+O_c(t^6).
\]
Thus the ratio $\delta(t)(1+t^2)/t^4$ tends to $\kappa>0$ as $t\to0$. It is strictly positive and continuous for $t\ne0$, since $\log(1+c\sin^2t)\le c\sin^2t<ct^2$ for $t\ne0$. Finally $0\le\log(1+c\sin^2t)\le\log(1+c)$, so the same ratio tends to $c>0$ as $|t|\to\infty$. Its positive infimum gives \eqref{eq:deficitcoercivity}.
\end{proof}

\begin{theorem}[Critical asymptotics and the anisotropic intensity profile]\label{thm:quarticcritical}
Let $Q=Q^T>0$, $v\in\R^d\setminus\{0\}$, $a>1$, and suppose
\begin{equation}\label{eq:boundarycondition}
 c\,v^TQ^{-1}v=1.
\end{equation}
Let $D=Q-cvv^T$, choose a unit vector $e\in\ker D$, write $E=e^\perp$, and denote $D_E=D|_E>0$ and $\nu=v\cdot e\ne0$.
Since $D$ is symmetric and $De=0$, the subspace $E=e^\perp$
is invariant under $D$, and the restriction $D_E$ is positive definite.
Set
\[
 h_N(x)=F_N(v\cdot x,a)e^{-x^TQx/(2N)}.
\]
For $w\in E$ and $z\in\R$ define
\begin{equation}\label{eq:criticaldensity}
 \rho_N(w,z)=\left|h_N\big(\sqrt N\,w+N^{3/4}ze\big)\right|^2.
\end{equation}
Then
\begin{equation}\label{eq:criticalL1}
 \rho_N\longrightarrow\rho_\infty(w,z)
 :=e^{-w^TD_Ew-\kappa\nu^4 z^4}
 \quad\text{in }L^1(E\times\R).
\end{equation}
In particular,
\begin{equation}\label{eq:quarticmass}
 \boxed{\displaystyle
 \|h_N\|_2^2\sim
 \frac{\pi^{(d-1)/2}}{\sqrt{\det D_E}}
 \frac{\Gamma(1/4)}{2(\kappa\nu^4)^{1/4}}
 N^{d/2+1/4}.}
\end{equation}
Here $\det D_E$ is the determinant of the matrix of the restriction $D|_E$ in any orthonormal basis of $E$; its value is basis-independent. For $d=1$ the determinant on the zero-dimensional space $E$ equals $1$.
Equivalently the normalized intensity measures
\[
 \frac{\rho_N(w,z)\,dw\,dz}{\int_{E\times\R}\rho_N(w,z)\,dw\,dz}
\]
converge in total variation to the probability measure proportional to $\rho_\infty(w,z)\,dw\,dz$.
\end{theorem}
\begin{proof}
The rank-one factorization in the proof of Theorem~\ref{thm:globalcritical} implies that $D$ is nonnegative, of rank $d-1$, and has a one-dimensional kernel; also $De=0$ implies $Qe=c(v\cdot e)v$, hence $\nu\ne0$. In particular $\nu\ne0$ guarantees that the quartic deficit controls the longitudinal null direction rather than vanishing there. The substitution
\[
 x=\sqrt N\,(w+N^{1/4}ze),\qquad w\in E,
\]
has Jacobian $N^{d/2+1/4}$. Put $y=w+N^{1/4}ze$. The identity $Q=D+cvv^T$ and the fact that $De=0$ give
$-y^TQy+c(v\cdot y)^2=-y^TDy=-w^TD_Ew$.
The exact modulus identity yields
\begin{align}
 |h_N(\sqrt N y)|^2
 &=\exp\left[-y^TQy+N\log\left(1+c\sin^2\frac{v\cdot y}{\sqrt N}\right)\right]\notag\\
 &=\exp\left[-w^TD_Ew-N\delta\left(\frac{v\cdot y}{\sqrt N}\right)\right].\label{eq:criticalexact}
\end{align}
Write $s_N=\nu z+N^{-1/4}v\cdot w$. Since $(v\cdot y)/\sqrt N=N^{-1/4}s_N$, the expansion in Lemma~\ref{lem:coercivity} implies, for each fixed $(w,z)$,
\[
 N\delta(N^{-1/4}s_N)\longrightarrow\kappa\nu^4z^4.
\]
This proves pointwise convergence in \eqref{eq:criticalL1}. We next construct an integrable majorant independent of $N$.

Let $\phi(r)=r^4/(1+r^2)$. By \eqref{eq:deficitcoercivity}, for $N\ge1$,
\begin{equation}\label{eq:criticaldeltauniform}
 N\delta(N^{-1/4}s_N)
 \ge k_c\frac{s_N^4}{1+N^{-1/2}s_N^2}
 \ge k_c\phi(s_N).
\end{equation}
We claim that there exist $C_1,C_2>0$, independent of $N,w,z$, with
\begin{equation}\label{eq:coercivepair}
 w^TD_Ew+k_c\phi(s_N)
 \ge C_1|w|^2+C_2\phi(z).
\end{equation}
Indeed, choose $\lambda>0$ such that $w^TD_Ew\ge\lambda|w|^2$. The function $\phi$ is even and increasing in $|r|$, since for $r>0$, $\phi'(r)=2r^3(2+r^2)/(1+r^2)^2>0$. It satisfies $\phi(r)\le r^2$ and $\phi(r)\le r^4$. Moreover, for all real $r,s$,
\[
 \phi(r+s)\le 16\big(\phi(r)+s^2\big).
\]
To check this, use $\phi(t)\le\min(t^2,t^4)$, split into $|r|\ge|s|$ and $|r|<|s|$, and in the former case use $|r+s|\le2|r|$ together with $\phi(2r)\le16\phi(r)$; in the latter case use $\phi(r+s)\le(r+s)^2\le4s^2$. Applying this with
$\nu z=s_N-N^{-1/4}v\cdot w$ and using $N^{-1/2}\le1$ gives
\[
 \phi(\nu z)\le16\big(\phi(s_N)+|v|^2|w|^2\big).
\]
Since $\nu\ne0$, the ratio $\phi(z)/\phi(\nu z)$, defined for $z\ne0$, has the finite limit $|\nu|^{-4}$ as $z\to0$ and the finite limit $|\nu|^{-2}$ as $|z|\to\infty$. Continuity between the two limits yields $\phi(z)\le C_\nu\phi(\nu z)$ for a finite $C_\nu$. Therefore
$\phi(z)\le C\{\phi(s_N)+|w|^2\}$, with $C$ independent of $N$. Put $m=\min\{\lambda,k_c\}>0$ and choose $0<\varepsilon<m$.
Then the preceding comparison implies
\[
 \lambda|w|^2+k_c\phi(s_N)
 \ge \varepsilon |w|^2+(m-\varepsilon)
 (|w|^2+\phi(s_N))
 \ge \varepsilon |w|^2+\frac{m-\varepsilon}{C}\phi(z),
\]
which proves \eqref{eq:coercivepair} with explicit positive constants.

Combining \eqref{eq:criticalexact}--\eqref{eq:coercivepair},
\[
 0\le\rho_N(w,z)\le
 e^{-C_1|w|^2-C_2\phi(z)}.
\]
This is integrable on $E\times\R$: as $|z|\to\infty$, $\phi(z)\sim z^2$, while its value near zero is harmless. Dominated convergence establishes \eqref{eq:criticalL1}. The same convergence gives the normalization and total-variation claim.

Finally, integrating \eqref{eq:criticalexact} and using the Jacobian,
\[
 N^{-d/2-1/4}\|h_N\|_2^2
 \longrightarrow\int_Ee^{-w^TD_Ew}\,dw
 \int_\R e^{-\kappa\nu^4z^4}\,dz.
\]
The first integral is $\pi^{(d-1)/2}/\sqrt{\det D_E}$. Substitution $u=(\kappa\nu^4)^{1/4}z$ and the Gamma-function identity $\int_0^\infty e^{-u^4}du=\Gamma(1/4)/4$ give the second integral as $\Gamma(1/4)/[2(\kappa\nu^4)^{1/4}]$. This proves \eqref{eq:quarticmass}.
\end{proof}

\begin{proposition}[Exact scalar normal form and independent check of the critical constant]
\label{prop:scalarnormalform}
For $Q=Q^T>0$ and $v\ne0$, put
$\sigma^2=v^TQ^{-1}v$ and $c=a^2-1$. The exact mass of the widening
Gaussian family is
\begin{equation}\label{eq:scalaridentity}
 \|h_N\|_2^2
 =\frac{(\pi N)^{d/2}}{\sqrt{\det Q}\,\sqrt{\pi\sigma^2}}
 \int_\R e^{-y^2/\sigma^2}
 \left(1+c\sin^2\frac{y}{\sqrt N}\right)^Ndy.
\end{equation}
At $c\sigma^2=1$, this identity independently gives
\begin{equation}\label{eq:scalarcriticalconstant}
 \lim_{N\to\infty}N^{-d/2-1/4}\|h_N\|_2^2
 =\frac{\pi^{(d-1)/2}\sqrt c}{\sqrt{\det Q}}
 \frac{\Gamma(1/4)}{2\kappa^{1/4}}.
\end{equation}
Moreover, with the notation of Theorem~\ref{thm:quarticcritical},
\begin{equation}\label{eq:detidentity}
 \det Q=c\nu^2\det D_E,
\end{equation}
so that~\eqref{eq:scalarcriticalconstant} is exactly the
constant in~\eqref{eq:quarticmass}.
\end{proposition}
\begin{proof}
Under the normalized Gaussian density proportional to
$e^{-y^TQy}\,dy$ on $\R^d$, the scalar random variable $v\cdot y$
has density
$(\pi\sigma^2)^{-1/2}e^{-s^2/\sigma^2}\,ds$.
Indeed, its characteristic function is
$\exp(-\sigma^2t^2/4)$, as follows by completing the Gaussian
square. Thus, changing variables $x=\sqrt N\,y$ and using the
exact modulus of $F_N$ proves~\eqref{eq:scalaridentity}.

At $c\sigma^2=1$, the substitution $y=N^{1/4}z$ reduces
\eqref{eq:scalaridentity} to
\[
 \|h_N\|_2^2
 =\frac{\pi^{(d-1)/2}\sqrt c}{\sqrt{\det Q}}
 N^{d/2+1/4}\int_\R e^{-N\delta(N^{-1/4}z)}\,dz.
\]
The coercivity estimate~\eqref{eq:deficitcoercivity} yields
\[
 N\delta(N^{-1/4}z)
 \ge k_c\frac{z^4}{1+N^{-1/2}z^2}
 \ge k_c\frac{z^4}{1+z^2}
 \qquad(N\ge1),
\]
which supplies an integrable majorant independent of $N$.
Since $N\delta(N^{-1/4}z)\to\kappa z^4$ pointwise, dominated
convergence and the Gamma integral yield
\eqref{eq:scalarcriticalconstant}.

To verify~\eqref{eq:detidentity}, use the orthogonal decomposition
$\R^d=E\oplus\R e$ and write $v=v_E+\nu e$.
Since $D|_E=D_E$ and $De=0$, the identity $Q=D+cvv^T$ has
block form
\[
 Q=\begin{pmatrix}
 D_E+c v_Ev_E^T & c\nu v_E\\
 c\nu v_E^T & c\nu^2
 \end{pmatrix}.
\]
The Schur complement of the lower-right entry $c\nu^2>0$
is exactly $D_E$. Taking determinants gives
$\det Q=c\nu^2\det D_E$. This agrees with the constant in
Theorem~\ref{thm:quarticcritical}, without invoking the
multidimensional Jacobian argument.
\end{proof}

\begin{corollary}[Sharp critical STFT energy law]\label{cor:quarticSTFT}
Let $\psi\in L^2(\R^d)\setminus\{0\}$. With the angular-frequency convention of Section~2,
\[
 \boxed{\displaystyle
 \|V_\psi h_N\|_{L^2(\R^{2d})}^2\sim
 (2\pi)^d\|\psi\|_2^2\,
 \frac{\pi^{(d-1)/2}\Gamma(1/4)}
 {2\sqrt{\det D_E}\,(\kappa\nu^4)^{1/4}}
 N^{d/2+1/4}.}
\]
\end{corollary}
\begin{proof}
Apply Moyal's identity $\|V_\psi h_N\|_2^2=(2\pi)^d\|\psi\|_2^2\|h_N\|_2^2$ and Theorem~\ref{thm:quarticcritical}.
\end{proof}

\begin{remark}[Scope of the result]
The limit \eqref{eq:criticalL1} is an anisotropically rescaled \emph{intensity} profile, not a complex-amplitude limit: the phases can oscillate on the critical scale. The STFT corollary concerns total phase-space energy and does not assert a local rescaled STFT profile. At the boundary, the spatial width is $N^{3/4}$ along $e$ and $N^{1/2}$ in transverse directions.
\end{remark}

\section{Entire growth and critical widening-Gaussian profiles}
\begin{proposition}[Exact complex-frequency growth]\label{prop:entire}
Let $Q=Q^T>0$ and $g_Q(x)=e^{-x^TQx/2}$. For real $\lambda,u\in\R^d$ and $z\in\C^d$,
\begin{equation}\label{gaussSTFT}
 V_{g_Q}(M_\lambda g_Q)(u,z)
 =\frac{\pi^{d/2}}{\sqrt{\det Q}}\exp\left[-\frac14u^TQu-\frac14(z-\lambda)^TQ^{-1}(z-\lambda)-\frac i2(z-\lambda)\cdot u\right].
\end{equation}
Put $B_Q=\frac14\lambda_{\max}(Q^{-1})$. For every fixed real $u,\lambda$, the entire function in \eqref{gaussSTFT} belongs to $A_{2,B}(\C^d)$ for every $B>B_Q$ and belongs to no $A_{p,B}(\C^d)$ with $p<2$ and finite $B$. It does not belong to $A_{2,B}$ for $B<B_Q$.
\end{proposition}
\begin{proof}
Completing the square in \eqref{STFT} proves \eqref{gaussSTFT}. Writing $z=\xi+i\eta$, its logarithmic modulus has quadratic part $\frac14\eta^TQ^{-1}\eta-\frac14\xi^TQ^{-1}\xi$ and terms of degree at most one. This gives the upper bound for $B>B_Q$. Setting $\xi=0$ and taking $\eta$ along a unit eigenvector for the largest eigenvalue of $Q^{-1}$ gives the lower growth $B_Q|\eta|^2+O(|\eta|)$, proving the remaining claims.\end{proof}

\begin{proposition}[Membership at the exact Gaussian growth type]\label{prop:gaussendpoint}
Let $E_{\max}$ be the top eigenspace of $Q^{-1}$. For fixed
$u,\lambda\in\mathbb R^d$, the Gaussian STFT in
\eqref{gaussSTFT} belongs to $A_{2,B_Q}(\mathbb C^d)$ if and only if
$\operatorname{proj}_{E_{\max}}u=0$.
\end{proposition}
\begin{proof}
Write $z=\xi+i\eta$, $A=Q^{-1}$, and use \eqref{gaussSTFT}.
Apart from factors independent of $\eta$, the logarithm of the
weighted modulus is
\[
 -\eta^T(B_Q I-A/4)\eta+\tfrac12u\cdot\eta
 -\tfrac14(\xi-\lambda)^TA(\xi-\lambda).
\]
The nonnegative matrix $B_Q I-A/4$ has kernel $E_{\max}$.
If the projection of $u$ onto this kernel is nonzero, choose
$\eta$ along that projection with the favorable sign; the supremum
diverges. If the projection vanishes, the linear term lies in the
orthogonal complement of the kernel, where the quadratic form is
positive definite. Completing the square on that complement bounds
its supremum; the $\xi$-term is nonpositive. This proves the claim.
\end{proof}

\begin{remark}Proposition \ref{prop:entire} concerns membership of each Gaussian STFT kernel. It does \emph{not} assert that the signed supershift converges in the global fixed-step $A_{2,B}$ norm. That stronger conclusion requires uniform control of the entire signed sum and should not be inferred from local uniform convergence.\end{remark}

\begin{proposition}[Critical scaling and a matrix threshold]\label{prop:critical}
Let $v\in\R^d\setminus\{0\}$, $Q=Q^T>0$, and
\[h_N(x)=F_N(v\cdot x,a)\exp\left(-\frac{x^TQx}{2N}\right).\]
Then uniformly on compact sets of $y\in\R^d$,
\begin{equation}\label{critlimit}
 |h_N(\sqrt N\,y)|^2\longrightarrow
 \exp\big(-y^T[Q-(a^2-1)vv^T]y\big).
\end{equation}
The limiting profile decays in every nonzero direction precisely when $Q-(a^2-1)vv^T$ is positive definite, equivalently
\begin{equation}\label{matrixthreshold} (a^2-1)v^TQ^{-1}v<1.\end{equation}
At equality it has a flat direction; when the inequality reverses it grows along some direction.
\end{proposition}
\begin{proof}
For $\theta=\sqrt N\,v\cdot y$, the exact modulus identity $|F_N(\theta,a)|^2=(1+(a^2-1)\sin^2(\theta/N))^N$ and $\sin(\theta/N)=\sin((v\cdot y)/\sqrt N)$ give
\[N\log\left(1+(a^2-1)\sin^2((v\cdot y)/\sqrt N)\right)\to(a^2-1)(v\cdot y)^2\]
uniformly on compacts. Multiplication by $\exp(-y^TQy)$ proves \eqref{critlimit}. The matrix criterion follows from the rank-one update identity $Q-(a^2-1)vv^T>0\iff (a^2-1)v^TQ^{-1}v<1$.\end{proof}
\begin{remark}This proposition alone proves only compact-uniform convergence of the rescaled intensity. The earlier global dichotomy and quartic theorem supply the separate tail control needed for total mass and $L^1$ convergence. None of these results is a local moving-frame STFT concentration theorem.\end{remark}

\section{A fixed-step entire-growth theorem for the Gaussian STFT}
The passage from pointwise convergence of entire STFT kernels to a common
Banach step is not automatic: the coefficients of the supershift have total
variation $a^N$. The next theorem gives a genuine fixed-step statement for
Gaussian-localized directional superoscillations, without estimating the
summands separately by their total variation.
\begin{theorem}[Fixed-step Gaussian STFT convergence]\label{thm:fixedstep}
Let $Q>0$ be symmetric, $v\in\R^d\setminus\{0\}$, $a>1$,
$g_Q(x)=\exp(-x^TQx/2)$, and
$f_N(x)=F_N(v\cdot x,a)g_Q(x)$, $f_a(x)=e^{ia v\cdot x}g_Q(x)$.
For fixed $u\in\R^d$ define
$\Phi_N(u,z)=V_{g_Q}f_N(u,z)$ and $\Phi_a(u,z)=V_{g_Q}f_a(u,z)$,
where $z\in\C^d$ is the complexified frequency. If
\[
 B> B_Q:=\frac14\lambda_{\max}(Q^{-1}),
\]
then
\[
 \|\Phi_N(u,\cdot)-\Phi_a(u,\cdot)\|_{A_{2,B}(\C^d)}\longrightarrow0.
\]
The convergence is uniform for $u$ in any prescribed compact subset of
$\R^d$. No rate in this fixed-step norm is asserted. The limiting kernel has quadratic order exactly two and has no
membership in $A_{2,B}$ for $B<B_Q$.
\end{theorem}
\begin{proof}
Write $z=\xi+i\eta$. For complex frequency the modulus of the
integrand is bounded by
\[
 |F_N(v\cdot x,a)|\,
 e^{-\frac12x^TQx-\frac12(x-u)^TQ(x-u)+\eta\cdot x}.
\]
Set $c=a^2-1$. The exact identity
$|F_N(\theta,a)|^2=(1+c\sin^2(\theta/N))^N$ gives
$|F_N(\theta,a)|\le \exp(c\theta^2/(2N))$.
Consequently, if
$Q_N=Q-\frac c{2N}vv^T>0$, completing the square yields
\[
 |\Phi_N(u,\xi+i\eta)|
 \le C_{Q,u}\int_{\R^d}
 \exp\{-x^TQ_Nx+(Qu+\eta)\cdot x\}\,dx.
\]
Completing the square gives
\[
 |\Phi_N(u,\xi+i\eta)|
 \le C'_{N,Q,u}\exp\left\{
 \frac14(Qu+\eta)^TQ_N^{-1}(Qu+\eta)
 \right\}.
\]
(The harmless factor $\exp(-u^TQu/2)$ is included in the prefactor.)
Choose $B_Q<B_0<B_1<B$. Since $Q_N\to Q$, for all sufficiently
large $N$ the preceding bound, together with Young's inequality for the
linear terms in $\eta$, implies
\[
 \sup_{u\in K,\xi\in\R^d}
 |\Phi_N(u,\xi+i\eta)|\le C_K e^{B_1|\eta|^2},
\]
for each compact $K\subset\R^d$; the same is true for $\Phi_a$.
Here the prefactor is uniform for $N$ large because $Q_N^{-1}$ remains
bounded. The bound is independent of $\xi$.

By Lemma~\ref{lem:global}, $f_N(x)\to f_a(x)$ pointwise and, for $N$
large, the difference times the fixed Gaussian is dominated by an
integrable Gaussian. On every compact set of $(u,z)$ the resulting STFT
integrands possess a common integrable majorant. In fact, for compact $K\subset\mathbb R^d$ and
$L\subset\mathbb C^d$, taking the supremum over $(u,z)\in K\times L$
inside the integral gives
\[
 \sup_{(u,z)\in K\times L}|\Phi_N(u,z)-\Phi_a(u,z)|
 \le \int_{\mathbb R^d}|f_N(x)-f_a(x)|G_{K,L}(x)\,dx,
\]
where $G_{K,L}(x)=\sup_{u\in K,z\in L}
 |g_Q(x-u)e^{-iz\cdot x}|$ is bounded by
$C\exp(-c|x|^2+C|x|)$ and is independent of $N$.
The right-hand side tends to zero by dominated convergence, because
$|f_N-f_a|$ has a common Gaussian majorant for all sufficiently large
$N$. Thus convergence is uniform on compact parameter sets. Notice that this argument uses dominated convergence only after taking the supremum over the compact parameter set inside the integral; pointwise dominated convergence alone would not suffice.

To pass to the Banach norm, choose $R>0$. On $|z|\le R$ uniform
convergence applies. On $|z|>R$ the uniform growth bound gives
\[
 (|\Phi_N(u,z)|+|\Phi_a(u,z)|)e^{-B|z|^2}
 \le 2C_K e^{-(B-B_1)|z|^2}.
\]
Its supremum tends to zero with $R\to\infty$, independently of large
$N$ and $u\in K$. Combining the compact and exterior regions proves the
claim. The last assertion is the exact Gaussian-kernel calculation of
Proposition~\ref{prop:entire}; along an imaginary eigenvector
corresponding to $\lambda_{\max}(Q^{-1})$ its modulus grows like
$\exp(B_Q|\eta|^2+O(|\eta|))$, which rules out $B<B_Q$.
\end{proof}
\begin{remark}[Endpoint is genuinely different]
The theorem proves strong convergence for $B>B_Q$ and nonmembership of
the limit for $B<B_Q$. At $B=B_Q$, the linear exponential factor and the
chosen spatial location $u$ can affect membership; no universal endpoint
convergence assertion is made. Likewise, the theorem is qualitative in the
fixed step and does not imply an $O(N^{-1})$ rate there without a further
weighted remainder estimate.
\end{remark}

\section{Failure of convergence at the critical Gaussian entire type}\label{sec:endpointfailure}
The fixed-step theorem above requires a strict increase of the Gaussian type. This is not merely an artifact of the proof. Even at the origin in the time variable the critical norm of the difference diverges exponentially.
\begin{theorem}[Sharp endpoint obstruction]\label{thm:endpointobstruction}
Let $a>1$, $g(x)=e^{-x^2/2}$, and $f_N(x)=F_N(x,a)g(x)$, $f_a(x)=e^{iax}g(x)$. Put $\Phi_N(z)=V_gf_N(0,z)$ and $\Phi_a(z)=V_gf_a(0,z)$. Then $\Phi_N,\Phi_a\in A_{2,1/4}(\C)$ for every fixed $N$, but
\[
\|\Phi_N-\Phi_a\|_{2,1/4}\ge \sqrt\pi\left[e^{-1/4}\left(\frac{1+a}{2}\right)^N-e^{-a^2/4}\right].
\]
In particular, the critical $A_{2,1/4}$-norm difference tends to infinity, although $\Phi_N\to\Phi_a$ locally uniformly and in every $A_{2,B}$ with $B>1/4$.
\end{theorem}
\begin{proof}
Direct Gaussian integration gives
\[
 \Phi_N(z)=\sqrt\pi\sum_{j=0}^N c_{j,N}(a)e^{-(z-s_{j,N})^2/4},\qquad
 \Phi_a(z)=\sqrt\pi e^{-(z-a)^2/4}.
\]
For real $s$ and $z=x+iy$, the modulus of $e^{-(z-s)^2/4}$ is $e^{y^2/4-(x-s)^2/4}$, hence the displayed finite sums belong to $A_{2,1/4}$; the same holds for $\Phi_a$. On the imaginary axis,
\[
 e^{-y^2/4}\Phi_N(iy)/\sqrt\pi
 =\sum_{j=0}^N c_{j,N}(a)e^{-s_{j,N}^2/4}e^{is_{j,N}y/2}
 =:T_N(y).
\]
The frequencies $s_{j,N}/2=1/2-j/N$ are distinct. Orthogonality on the interval $[0,4\pi N]$ gives
\[
 \frac1{4\pi N}\int_0^{4\pi N}|T_N(y)|^2\,dy
 =\sum_{j=0}^N |c_{j,N}(a)|^2e^{-s_{j,N}^2/2}.
\]
Therefore
\[
 \sup_y |T_N(y)|\ge |c_{0,N}(a)|e^{-1/4}
 =e^{-1/4}\left(\frac{1+a}{2}\right)^N.
\]
On the other hand,
$e^{-y^2/4}|\Phi_a(iy)|=\sqrt\pi e^{-a^2/4}$ for all real $y$. The reverse triangle inequality, followed by the definition of the weighted supremum norm, proves the assertion.
\end{proof}
\begin{remark}
The obstruction concerns the \emph{fixed Banach step at the critical type}, not the larger space $A_2$ or any step with $B>1/4$. In particular, individual membership at critical type does not imply convergence there. In higher dimensions, restriction to an imaginary-frequency line in the maximal eigenspace permits the same orthogonality argument if the projection of $v$ onto that eigenspace is nonzero. The result is not claimed for arbitrary anisotropic matrices or translations, and no general assertion is made when that projection vanishes.
\end{remark}

\section{Compact windows and modulation-space convergence}
Fix $v\in\R^d\setminus\{0\}$ and $g\in C_c^\infty(\R^d)$.
Set $w_N(x)=F_N(v\cdot x,a)g(x)$,
$w_a(x)=e^{ia v\cdot x}g(x)$, and
\[
 v_a(x)=\frac{a^2-1}{2}(v\cdot x)^2w_a(x).
\]
This directional model is sufficient for a genuinely multidimensional
Gabor-frame perturbation theorem. For a coupled nonlinear frequency curve
$\Lambda(s)=(s^{p_1},\dots,s^{p_d})$, an analytic supershift theorem
with uniform derivative estimates is needed; it is not silently inferred
from pointwise convergence here.

\begin{proposition}[Full smooth-window expansion]\label{prop:compact}
For every nonnegative integer $k$ and every $m\ge1$,
\[
 \left\|w_N-e^{ia v\cdot x}g
 \sum_{j=0}^m N^{-j}P_j(v\cdot x)\right\|_{C^k(\R^d)}
 \le C_{k,m}N^{-m-1}.
\]
In particular,
\begin{equation}\label{M1exp}
 \|w_N-w_a-N^{-1}v_a\|_{M^1}\le CN^{-2}.
\end{equation}
\end{proposition}
\begin{proof}
Let $K=\supp g$. For $x\in K$, $|v\cdot x|$ is uniformly bounded.
The analytic Taylor expansion used in Lemma~\ref{lem:global} therefore
has a remainder bounded by $C_{K,m}N^{-m-1}$, uniformly in $x\in K$.
For $\beta\in\N_0^d$, differentiation of
$F_N(v\cdot x,a)=\sum_jc_{j,N}e^{is_{j,N}v\cdot x}$
can equivalently be obtained by differentiating the analytic function
$F_N(\theta,a)$ in $\theta$: Cauchy's estimate on a fixed complex
neighborhood of $\{v\cdot x:x\in K\}$ bounds the derivatives of the
Taylor remainder by $C_{\beta,K,m}N^{-m-1}$.
Leibniz's rule and the bounded derivatives of $g$ give the stated
$C^k$ estimate, including at the boundary of $K$.

Here is a direct justification of the $M^1$ bound. Fix a Schwartz
window $\phi$ and an integer $L>2d$. If $r$ is supported in $K$,
repeated integration by parts in the defining STFT yields
\[
 |V_\phi r(u,\omega)|
 \le C_L(1+|\omega|)^{-L}
 \sum_{|\beta|\le L}\int_K
 \left|\partial_x^\beta\{r(x)\overline{\phi(x-u)}\}\right|dx.
\]
Integrating this expression first in $\omega$ and then in $u$ and
using the integrability of every Schwartz derivative of $\phi$ gives
\[
 \|r\|_{M^1}\le C_{K,L,\phi}\|r\|_{C^L}.
\]
Apply the inequality to $r=w_N-w_a-N^{-1}v_a$ and the already proved
$C^L$ estimate.\end{proof}

\begin{remark}[Nonlinear frequency curves]
For $\Lambda(s)=(s^{p_1},\ldots,s^{p_d})$, the formally expected first
correction is
\[
 \frac{1-a^2}{2}\bigl(i\Lambda''(a)\cdot x
 -(\Lambda'(a)\cdot x)^2\bigr)e^{i\Lambda(a)\cdot x}g(x).
\]
A rigorous $O(N^{-2})$ remainder for arbitrary such curves requires
uniform analytic supershift estimates on the compact support. This
possible extension is distinguished from the proved directional result.
\end{remark}

\section{Gabor frames: valid bounds and first variations}
Let $\Gamma$ be a fixed full-rank time--frequency lattice and suppose $\mathcal G(w_a,\Gamma)$ is a frame with bounds $A,B>0$. We work with the compactly supported windows from Proposition \ref{prop:compact}; the assertion is conditional on the actual frame property, not on nonvanishing of the underlying multiplier. In particular, the later indicator-window example is not an instance of the $M^1$ hypothesis here, and no Gaussian-window frame property is inferred from this theorem. Proposition~\ref{prop:smoothframe} below verifies the hypothesis for an explicit smooth compactly supported window and an oversampled lattice.
\begin{theorem}[Frame perturbation and reconstruction asymptotics]\label{thm:frames}
Let $w_N,w_a,v_a$ be the windows in Proposition~\ref{prop:compact},
and set $\rho_N:=w_N-w_a-N^{-1}v_a$. Suppose
$\|\rho_N\|_{M^1}\le CN^{-2}$, as established there.
If $K_\Gamma$ is the operator bound in~\eqref{M1frame}, one may take
$K=2K_\Gamma(\|v_a\|_{M^1}+C)$ for sufficiently large $N$; for all sufficiently large $N$, $\mathcal G(w_N,\Gamma)$ is a Gabor frame with bounds
\begin{equation}\label{framebounds}
 (\sqrt A-K/N)^2\le A_N,\qquad B_N\le(\sqrt B+K/N)^2.
\end{equation}
Let $S_N=D_{w_N}C_{w_N}$, $S_a=D_{w_a}C_{w_a}$, and
\[\mathcal K_a=D_{v_a}C_{w_a}+D_{w_a}C_{v_a}.\]
Then in $\mathcal B(L^2(\R^d))$,
\begin{align}
 S_N&=S_a+N^{-1}\mathcal K_a+O(N^{-2}),\label{Sexp}\\
 S_N^{-1}&=S_a^{-1}-N^{-1}S_a^{-1}\mathcal K_a S_a^{-1}+O(N^{-2}).\label{Sinverse}
\end{align}
The canonical dual windows satisfy in $L^2$,
\begin{equation}\label{dual}
 S_N^{-1}w_N=S_a^{-1}w_a+\frac1N\big(S_a^{-1}v_a-S_a^{-1}\mathcal K_a S_a^{-1}w_a\big)+O_{L^2}(N^{-2}).
\end{equation}
\end{theorem}
\begin{proof}
Write $w_N=w_a+N^{-1}v_a+\rho_N$, with precisely the
compactly supported windows $w_N,w_a,v_a$ defined in
Proposition~\ref{prop:compact}, where
$\|\rho_N\|_{M^1}\le C N^{-2}$. The dependence of the analysis
operator on the window is conjugate-linear under the convention
$C_wf=(\langle f,\pi(\gamma)w\rangle)_\gamma$ with
$\langle f,g\rangle=\int f\overline g$; for the real scalar
$N^{-1}$ the expansion is unaffected. More precisely,
\[
 C_{w_N}=C_{w_a}+N^{-1}C_{v_a}+C_{\rho_N},\qquad
 D_{w_N}=D_{w_a}+N^{-1}D_{v_a}+D_{\rho_N}.
\]
The bounds in~\eqref{M1frame} imply
$\|C_{\rho_N}\|+\|D_{\rho_N}\|\le 2K_\Gamma C N^{-2}$.
In particular, if $\delta_N=\|C_{w_N}-C_{w_a}\|$, then
$\delta_N\le K/N$. For every $f\in L^2$,
\[
 \|C_{w_N}f\|_{\ell^2}
 \ge \|C_{w_a}f\|_{\ell^2}-\delta_N\|f\|_2
 \ge(\sqrt A-\delta_N)\|f\|_2,
\]
whereas the triangle inequality gives
$\|C_{w_N}f\|_{\ell^2}\le(\sqrt B+\delta_N)\|f\|_2$.
For $N$ large enough that $K/N<\sqrt A$, this proves the
frame inequalities and the claimed bounds.

Expand the product using the preceding operator identities:
\[
 \begin{split}
 S_N={}&S_a+N^{-1}(D_{v_a}C_{w_a}+D_{w_a}C_{v_a})\\
 &+N^{-2}D_{v_a}C_{v_a}
 +D_{\rho_N}C_{w_a}+D_{w_a}C_{\rho_N}
 +O_{\mathcal B(L^2)}(N^{-3}).
 \end{split}
\]
Each displayed term after the first-order correction is
$O(N^{-2})$ in operator norm, proving~\eqref{Sexp}.
Let $E_N=S_N-S_a$ and $T_N=S_a^{-1}E_N$.
Since $\|T_N\|=O(N^{-1})$, eventually $\|T_N\|<1/2$.
The identity
\[
 (I+T_N)^{-1}=I-T_N+T_N^2(I+T_N)^{-1}
\]
therefore implies
\[
 S_N^{-1}=S_a^{-1}-S_a^{-1}E_NS_a^{-1}+O(N^{-2}),
\]
with all errors in operator norm. Substituting
$E_N=N^{-1}\mathcal K_a+O(N^{-2})$ gives~\eqref{Sinverse}.
Finally multiply this operator expansion by
$w_N=w_a+N^{-1}v_a+O_{L^2}(N^{-2})$.
All cross terms containing two first-order factors are $O(N^{-2})$,
which yields~\eqref{dual}.\end{proof}

\begin{proposition}[An explicit smooth window realizing the frame hypothesis]
\label{prop:smoothframe}
In dimension one, choose $\rho\in C_c^\infty(\R)$ such that
$\rho\ge0$, $\operatorname{supp}\rho\subset[-1,1]$, and
$\rho(x)>0$ for $|x|<1$; for instance,
\[
 \rho(x)=\begin{cases}
 \exp\bigl(-1/(1-x^2)\bigr),&|x|<1,\\
 0,&|x|\ge1.
 \end{cases}
\]
Define
\[
 R(x)=\sum_{k\in\Z}\rho(x-k)^2,
 \qquad g(x)=\rho(x)R(x)^{-1/2}.
\]
Then $g\in C_c^\infty(\R)$ and
$\sum_{k\in\Z}|g(x-k)|^2=1$ for all $x\in\R$.
For any $a>1$ and $v\ne0$, the Gabor system
\[
 \mathcal G(w_a,\Gamma)=
 \{M_{\pi m}T_kw_a:k,m\in\Z\},
 \quad w_a(x)=e^{iavx}g(x),\quad
 \Gamma=\Z\times\pi\Z,
\]
is a tight frame for $L^2(\R)$ with frame bound $2$.
The superoscillatory windows $w_N(x)=F_N(vx,a)g(x)$ belong to
$M^1(\R)$, and the preceding frame-perturbation theorem applies to
this explicit smooth-window family. In fact, every $w_N$ generates
a frame with frame operator
\[
 (S_Nf)(x)=2\sum_{k\in\Z}
 |F_N(v(x-k),a)|^2|g(x-k)|^2 f(x),
\]
and lower frame bound at least $2$. For all sufficiently large $N$,
its upper frame bound can be taken as $2(1+C_{a,v}/N)$.
\end{proposition}
\begin{proof}
The sum defining $R$ is locally finite and smooth, is $1$-periodic,
and is strictly positive: for each $x$ there is an integer $k$ with
$|x-k|\le1/2$. Thus $R^{-1/2}$ is smooth and bounded, and $g$ has the
asserted support and regularity. Periodicity gives
\[
 \sum_k|g(x-k)|^2=\frac{\sum_k\rho(x-k)^2}{R(x)}=1.
\]
For the translated window $q=T_kw_a$, its support is contained in
$[k-1,k+1]$. We apply Parseval on this interval of length $2$ to
$f\overline q$, extended by zero outside $[k-1,k+1]$; the functions
$\{2^{-1/2}e^{i\pi mx}\}_{m\in\Z}$ form an orthonormal basis
there (the translation by the integer $k$ changes them only by a
unimodular factor). More generally, for any function $q$ supported
in an interval of length $2$, Parseval gives, for $f\in L^2(\R)$,
\[
 \sum_{m\in\Z}|\langle f,M_{\pi m}q\rangle|^2
 =2\int_\R|f(x)|^2|q(x)|^2\,dx.
\]
Apply this to $q=T_kw_a$ and sum in $k$; all terms are nonnegative,
so Tonelli's theorem yields
\[
 \sum_{k,m\in\Z}|\langle f,M_{\pi m}T_kw_a\rangle|^2
 =2\int_\R|f(x)|^2\sum_k|g(x-k)|^2\,dx
 =2\|f\|_2^2.
\]
(The modulation $e^{iavx}$ does not change the squared modulus.)
The same calculation with $w_N$ proves the displayed multiplication
formula for $S_N$. Since for real $\theta$,
$|F_N(\theta,a)|^2=(1+(a^2-1)\sin^2(\theta/N))^N\ge1$,
the lower bound is at least $2$. On $\operatorname{supp}g$ we have
$|vx|\le|v|$, and the bound
\[
 |F_N(vx,a)|^2\le e^{(a^2-1)|v|^2/N}
 \le1+C_{a,v}/N
\]
holds for all large $N$. Summing the squared window translates gives
the claimed upper bound. Finally, Proposition~\ref{prop:compact}
provides the expansion in $M^1$, so this is also an explicit instance
of Theorem~\ref{thm:frames}. This is the smooth compact-support
(``painless'') frame mechanism; see \cite{DGM}.
\end{proof}

\begin{proposition}[Failure of a universal Zak lower bound]\label{prop:zak}
For any $a>1$ and any $N$, there exists a nonzero $g\in C_c^\infty(\R)$ such that the unit-lattice Zak transform of $F_N(\cdot,a)g$ vanishes on a set of positive measure. Hence there is no positive uniform Zak lower bound for arbitrary nonzero $L^2$ windows.
\end{proposition}
\begin{proof}
Take $0\ne g\in C_c^\infty((0,1/4))$. With $Zf(u,\eta)=\sum_{k\in\Z}f(u-k)e^{ik\eta}$, $Z(F_Ng)(u,\eta)=0$ for all $u\in(1/2,3/4)$ and all $\eta$.\end{proof}
\begin{remark}The Zak-transform proposition is a correction to the universal lower-bound assertion appearing in discussions of superoscillatory Gabor windows; it does not imply that all of the existing frame results are false. The conclusion $\inf|Zf|>0$ cannot follow merely from $g\ne0$, nor can cancellation in a signed sum be excluded term by term. We do not use that argument.\end{remark}
\begin{remark}[Schatten classes] If $S$ is a frame operator on infinite-dimensional $L^2$ with $S\ge AI$, $A>0$, then $S$ cannot be compact. In particular $S\notin\mathcal S_p$ for any finite $p$. Any Schatten assertion must concern a genuinely compact difference or localization, not the invertible frame operator itself.\end{remark}

\section{A transparent critical-density frame model}
In contrast to the false statement that every nonzero window generates a nonvanishing Zak transform, a concrete window permits exact frame bounds. This model also shows that a lower frame bound need not deteriorate when the upper one depends on $N$.
\begin{proposition}[Exact frame bounds for a painless superoscillatory window]\label{prop:painless}
Fix $a>1$ and set
\[
 w_N(x)=F_N(x,a)\mathbf1_{[0,1)}(x).
\]
With $M_{2\pi m}f(x)=e^{2\pi imx}f(x)$, the Gabor system
\[
 \mathcal G(w_N)=\{M_{2\pi m}T_kw_N:k,m\in\Z\}
\]
is a frame for $L^2(\R)$. Its optimal lower bound is $1$, and its optimal upper bound is
\[
 B_N=\operatorname*{ess\,sup}_{0\le x<1}|F_N(x,a)|^2\le a^{2N}.
\]
Writing $w_a(x)=e^{iax}\mathbf1_{[0,1)}(x)$, its frame operator is the identity, and
\[
 \|S_{w_N}-I\|_{L^2\to L^2}
 =\operatorname*{ess\,sup}_{0\le x<1}(|F_N(x,a)|^2-1)=O_a(N^{-1}).
\]
Moreover the canonical dual window is given exactly by
\[
 \widetilde w_N(x)=\frac{w_N(x)}{|F_N(x,a)|^2}\quad(0\le x<1),
 \qquad \widetilde w_N(x)=0\quad(x\notin[0,1)),
\]
and $\|\widetilde w_N-w_a\|_2=O_a(N^{-1})$.
\end{proposition}
\begin{proof}
The translates of $[0,1)$ partition the real line up to null sets.
For each $k\in\Z$, Parseval's identity for the exponentials
$\{e^{2\pi imx}:m\in\Z\}$ on $[k,k+1)$ gives
\[
 \sum_{m\in\Z}|\langle f,M_{2\pi m}T_kw_N\rangle|^2
 =\int_k^{k+1}|f(x)|^2|F_N(x-k,a)|^2\,dx.
\]
Summing in $k$ shows that the frame operator acts on $[k,k+1)$
as multiplication by $|F_N(x-k,a)|^2$. The exact modulus identity
\[
 |F_N(s,a)|^2=\left(1+(a^2-1)\sin^2(s/N)\right)^N
\]
gives the stated bounds and shows that the essential infimum is $1$
(the multiplier approaches $1$ as $s\downarrow0$ on sets of positive measure).
It also gives the exact inverse multiplier, hence the dual-window formula.

For $0\le s\le1$ we have
\[
 0\le|F_N(s,a)|^2-1
 \le e^{(a^2-1)/N}-1=O_a(N^{-1}),
\]
using $\sin(s/N)\le1/N$. This proves operator-norm convergence.
Uniformly on $[0,1]$, analytic expansion near the origin gives
$F_N(s,a)=e^{ias}+O_a(N^{-1})$. Since $|F_N(s,a)|\ge1$,
\[
 \left|\frac{F_N(s,a)}{|F_N(s,a)|^2}-e^{ias}\right|
 \le|F_N(s,a)-e^{ias}|+||F_N(s,a)|^{-2}-1|
 =O_a(N^{-1}),
\]
which gives the asserted $L^2$ estimate.
\end{proof}

\begin{remark}
This critical-density example is elementary; its role is to exhibit exact reconstruction under explicit support assumptions. It is not a general replacement for the perturbative frame theorem, and the window $\mathbf1_{[0,1)}$ is not in $M^1$.
\end{remark}

\section{Higher-order Schr\"odinger evolution and uniform electric fields}
Let $P:\R^d\to\R$ be a real polynomial (not necessarily separable or even) and $b\in\R^d$. The operator
\[H_{P,b}=P(-i\nabla)-b\cdot x\]
has a self-adjoint realization through its explicit Fourier-side unitary propagator. We write $U_{P,b}(t)$ for the associated unitary group.
For real polynomial $P$, the Fourier-side formula below is strongly
continuous in $t$: translations are strongly continuous on $L^2$,
and the phase factors converge pointwise and have modulus one, so
dominated convergence applies. The characteristic integral obeys the
cocycle identity, yielding the group law. Stone's theorem therefore gives a self-adjoint generator.
The identification of this generator with the differential expression,
and the Schwartz-core property, are proved in the following proposition.

\begin{proposition}[Schwartz core for polynomial electric-field Hamiltonians]
\label{prop:electriccore}
Let $P$ be a real polynomial on $\R^d$ and let $b\in\R^d$.
The operator $P(D)-b\cdot x$ on $\mathcal S(\R^d)$ is essentially
self-adjoint. Its closure is the self-adjoint generator of the
strongly continuous unitary group defined by the Fourier
characteristic formula~\eqref{fourierprop}.
\end{proposition}
\begin{proof}
With the Fourier convention of Section~2, the differential
expression becomes, on $\mathcal S(\R^d)$,
\[
 \widehat H=P(\xi)-i b\cdot\nabla_\xi.
\]
When $b=0$, this is multiplication by the real polynomial $P$.
The space $C_c^\infty(\R^d)$ is a core for the maximal real
multiplication operator: cutoffs approximate any function in its
graph norm $\|f\|_2+\|Pf\|_2$, and compactly supported functions
can then be mollified. Hence $\mathcal S$ is also a core.

Suppose $b\ne0$. Put
\[
 t(\xi)=\frac{b\cdot\xi}{|b|^2},\qquad
 \xi_\perp=\xi-t(\xi)b,\qquad
 \Phi(\xi)=\int_0^{t(\xi)}P(\xi_\perp+s b)\,ds.
\]
The function $\Phi$ is a real polynomial because $P$ is a
polynomial and $t,\xi_\perp$ are linear in $\xi$.
Translation of $\xi$ by $s b$ changes $t$ by $s$ while leaving
$\xi_\perp$ fixed. Consequently
\[
 b\cdot\nabla_\xi\Phi(\xi)=P(\xi).
\]
Multiplication by $e^{i\Phi}$ is unitary on $L^2(\R^d)$ and maps
$\mathcal S(\R^d)$ bijectively onto itself, since its derivatives
are polynomial multiples of $e^{i\Phi}$. Direct differentiation gives
\[
 e^{-i\Phi}(-i b\cdot\nabla_\xi)e^{i\Phi}
 =-i b\cdot\nabla_\xi+P(\xi)
 \qquad\text{on }\mathcal S(\R^d).
\]
The constant-coefficient transport operator
$-i b\cdot\nabla_\xi$ is essentially self-adjoint on
$\mathcal S(\R^d)$, by the Fourier transform or the standard
self-adjointness of a real linear multiplication operator.
Unitary equivalence, including invariance of $\mathcal S$ under
the gauge factors, proves essential self-adjointness of $\widehat H$
on $\mathcal S$.

The Fourier characteristic formula defines a strongly continuous
unitary group: translations are strongly continuous, the multiplier
has modulus one, and dominated convergence applies to its phase as
$t\to0$. The characteristic cocycle gives the group law. On
$\mathcal S$, its derivative at $t=0$ is $-i\widehat H$.
By Stone's theorem and essential self-adjointness, its generator
is exactly the closure of $\widehat H$, as claimed.
\end{proof}

\begin{theorem}[Exact STFT formula and quantitative persistence]\label{thm:electric}
For $f\in\Sspace(\R^d)$, the solution of $i\partial_t\psi=H_{P,b}\psi$, $\psi(0)=f$, has Fourier transform
\begin{equation}\label{fourierprop}
 \widehat\psi(t,\xi)=\exp\left[-i\int_0^tP(\xi-tb+sb)\,ds\right]\widehat f(\xi-tb).
\end{equation}
Put $\Theta_t(\xi)=\int_0^tP(\xi+sb)\,ds$ and let $A_t$ be the Fourier multiplier $\widehat{A_tf}(\xi)=e^{-i\Theta_t(\xi)}\widehat f(\xi)$. Then
\begin{equation}\label{factor}
 U_{P,b}(t)f=M_{tb}A_tf,\qquad
 V_h(U_{P,b}(t)f)(u,\omega)=V_h(A_tf)(u,\omega-tb).
\end{equation}
If $f_N,f_a,r_1,r_2$ are as in Theorem \ref{thm:STFT}, then, for every $2\le q\le\infty$ and $h\in L^2$,
\begin{multline}\label{evolve}
 \sup_{t\in\R}\big\|V_hU_{P,b}(t)f_N-V_hU_{P,b}(t)f_a
 -N^{-1}V_hU_{P,b}(t)r_1\\-N^{-2}V_hU_{P,b}(t)r_2\big\|_{L^q}=O(N^{-3}).
\end{multline}
For the more general case of a real measurable $P$, suppose the
oriented characteristic integral
$\int_0^tP(\eta+sb)\,ds$ is finite and measurable for almost
every $\eta$ and every fixed $t$, with sufficient local integrability
along almost every characteristic to define the phase consistently.
Under these assumptions, for each fixed $t$ the Fourier formula
defines an isometry (indeed a unitary operator) on $L^2$ and gives
the stated unweighted STFT bound. We do not assert a strongly continuous
one-parameter group for this measurable-symbol extension without
separately verifying the cocycle identity and strong continuity.
The polynomial hypothesis is retained for the Schwartz-space and
weighted estimates below.
\end{theorem}
\begin{proof}
Using the angular-frequency Fourier convention,
$\widehat{x_j\psi}=i\partial_{\xi_j}\widehat\psi$,
so the Fourier transform of the equation is
\[
 \partial_t\widehat\psi(t,\xi)
 +b\cdot\nabla_\xi\widehat\psi(t,\xi)
 =-iP(\xi)\widehat\psi(t,\xi).
\]
Fix $\eta\in\R^d$ and follow the characteristic
$\xi(s)=\eta+sb$. Then
\[
 \frac{d}{ds}\widehat\psi(s,\eta+sb)
 =-iP(\eta+sb)\widehat\psi(s,\eta+sb).
\]
Integrating this scalar equation gives
\[
 \widehat\psi(t,\eta+tb)=
 \exp\left(-i\int_0^t P(\eta+sb)\,ds\right)\widehat f(\eta).
\]
Substitution $\eta=\xi-tb$ proves~\eqref{fourierprop}, also
for negative $t$ with the usual oriented integral. Since $P$ is
real valued, the exponential factor has modulus one. Therefore the
propagator is an isometry on the Fourier side by Plancherel and
translation invariance of Lebesgue measure; its inverse is obtained
by reversing time. This constructs the unitary evolution on $L^2$.

By definition, $\widehat{A_tf}(\eta)=e^{-i\Theta_t(\eta)}\widehat f(\eta)$.
Since $\widehat{M_{tb}A_tf}(\xi)=\widehat{A_tf}(\xi-tb)$,
we obtain $U_{P,b}(t)=M_{tb}A_t$. Direct substitution in the
STFT integral gives
\[
 V_h(M_{tb}k)(u,\omega)
 =\int k(x)\overline{h(x-u)}e^{-i(\omega-tb)\cdot x}\,dx
 =V_hk(u,\omega-tb),
\]
proving~\eqref{factor}.

Let $E_N=f_N-f_a-N^{-1}r_1-N^{-2}r_2$.
Lemma~\ref{lem:global} yields $\|E_N\|_2\le CN^{-3}$.
By linearity, the expression inside the norm in~\eqref{evolve}
is $V_h(U_{P,b}(t)E_N)$. Thus for every $t\in\R$,
\[
 \|V_h(U_{P,b}(t)E_N)\|_q
 \le(2\pi)^{d/q}\|h\|_2\|U_{P,b}(t)E_N\|_2
 =(2\pi)^{d/q}\|h\|_2\|E_N\|_2.
\]
The right-hand side is independent of $t$, completing the proof.
For the measurable-symbol extension in the statement, the same
unit-modulus phase and translation argument gives the fixed-time
$L^2$ and STFT bounds under the explicit characteristic-integrability
assumptions; no Schwartz-space or unverified strongly continuous
group assertion is made for such a general symbol.\end{proof}
\begin{remark}[Phase-space dynamics] For the Hamiltonian $\mathcal H(x,\xi)=P(\xi)-b\cdot x$, classical trajectories satisfy $\xi(t)=\xi_0+tb$ and $x(t)=x_0+\int_0^t\nabla P(\xi_0+sb)\,ds$. These trajectories motivate weighted or moving-window STFT estimates. Such genuinely microlocal bounds are not claimed by Theorem \ref{thm:electric}, whose norm estimate follows from unitarity.\end{remark}
\begin{remark}[Infinite-order expressions] If $P(\xi)=\sum_\alpha a_\alpha\xi^\alpha$ is merely analytic near a bounded frequency set, its derivative series need not act termwise on Gaussian signals (which are not bandlimited). Fourier-multiplier functional calculus remains meaningful when $P$ is real and measurable on all of $\R^d$. Identifying it with a convergent infinite differential series requires an additional domain theorem. Classical AU-space results for analytic supershifts are treated in \cite{AMS}; they are not new consequences of the present $L^2$ argument.\end{remark}
\begin{remark}[The original electric-field data] For $0\le\tau<1$, a window $g(x_0-(1-\tau)x)$ is in $L^2$ whenever $g\in L^2$, and fits this theorem in the Gaussian case after a change of centre and covariance. At $\tau=1$ it is constant and produces plane-wave data outside $L^2$, so a distributional or analytic-growth formulation is necessary.\end{remark}

\section{Polynomially weighted STFT persistence in an electric field}
The unweighted estimate in Theorem~\ref{thm:electric} follows from
unitarity. Polynomial phase-space weights require a stronger topology and
regularity of the spectral symbol. We give a direct finite-time estimate
which is not restricted to quadratic dispersion.
\begin{definition}
We use throughout the angular-frequency convention
\[
 V_hf(u,\omega)=\int_{\R^d} f(x)\overline{h(x-u)}
 e^{-i\omega\cdot x}\,dx,
 \qquad \|V_hf\|_{L^2(\R^{2d})}^2
 =(2\pi)^d\|h\|_2^2\|f\|_2^2.
\]
For $m\in\N_0$ and a Schwartz window $h$ define
\[
 \|f\|_{\mathcal V_m(h)}
 :=\left(\int_{\R^{2d}}(1+|u|+|\omega|)^{2m}
 |V_hf(u,\omega)|^2\,du\,d\omega\right)^{1/2}.
\]
\end{definition}
\begin{lemma}[Schwartz control of weighted STFTs]\label{lem:weighted}
For every integer $m\ge0$ and fixed $h\in\mathcal S(\R^d)$ there
exist an integer $M$ and a constant $C$ such that
\[
 \|f\|_{\mathcal V_m(h)}\le C
 \sum_{|\alpha|+|\beta|\le M}\|x^\alpha\partial^\beta f\|_2,
 \qquad f\in\mathcal S(\R^d).
\]
\end{lemma}
\begin{proof}
The inequality $(1+|u|+|\omega|)^m\le
C_m\sum_{|\alpha|+|\beta|\le m}|u^\alpha\omega^\beta|$
reduces matters to monomials. In the defining STFT integral,
$u_j=x_j-(x_j-u_j)$, so multiplying by $u^\alpha$ expresses the
STFT as a finite linear combination of STFTs of $x^\gamma f$
with windows $y^\delta h(y)$. Integration by parts in $x$ gives
\[
 \omega_j V_kf(u,\omega)
 =-i\,V_k(\partial_j f)(u,\omega)
   -i\,V_{\partial_jk}f(u,\omega),
\]
with the convention that differentiating the conjugated window gives
its conjugated derivative. Iterating these two identities and applying
Moyal's formula to every term yields the result, with $M=2m$ sufficient.
\end{proof}
\begin{lemma}[All-order Gaussian remainder bounds]\label{lem:schwartzremainder}
For the directional Gaussian signals of Lemma~\ref{lem:global}, put
\[
E_N(x)=g(x)\left\{F_N(v\cdot x,a)-e^{ia v\cdot x}
\big[1+N^{-1}P_1(v\cdot x)+N^{-2}P_2(v\cdot x)\big]\right\}.
\]
For every pair of multi-indices $\alpha,\beta\in\mathbb N_0^d$,
there exist constants $C_{\alpha,\beta}>0$ and
$N_{\alpha,\beta}\in\mathbb N$ (depending on
$a,v,Q,x_0,\alpha,\beta$ but not on $N$) such that
\[
 \|x^\alpha\partial_x^\beta E_N\|_{L^2(\mathbb R^d)}
 \le C_{\alpha,\beta}N^{-3},
 \qquad N\ge N_{\alpha,\beta}.
\]
\end{lemma}
\begin{proof}
Write $\theta=v\cdot x$ and $c=a^2-1$. Let
\[
 L(z)=\log(\cos z+ia\sin z),\qquad L(0)=0,
\]
which is holomorphic on $|z|<r_0$ for some $r_0>0$.
For $t$ near zero put
\[
 H(t,\zeta)=\exp\left(\frac{L(t\zeta)-ia t\zeta}{t}\right),
 \qquad H(0,\zeta)=1.
\]
For every fixed $\theta$, this defines a jointly holomorphic function of $(t,\zeta)$ on a neighborhood of the polydisk used below. Indeed, if $L(z)=ia z+\sum_{j\ge2}\ell_j z^j$ on $|z|<r_0$, then
\[
 \frac{L(t\zeta)-ia t\zeta}{t}
 =\sum_{j\ge2}\ell_j t^{j-1}\zeta^j
\]
converges normally on sufficiently small polydisks, explicitly removing the singularity at $t=0$. Since
$L(z)-iaz=O(z^2)$, there is $C>0$ such that
\[
 \left|\frac{L(t\zeta)-ia t\zeta}{t}\right|
 \le C|t||\zeta|^2\qquad (|t\zeta|\le r_0/2).
\]
Fix $0<\varepsilon<1/2$ and a small constant $\eta>0$ such that $2\eta<r_0$. For real $\theta$ and
$|\zeta-\theta|\le1$, put
$r_\theta=\eta(1+|\theta|)^{-2}$.
Because $|\zeta|\le1+|\theta|$ on this disk, we have $|t|\,|\zeta|^2\le\eta$ for $|t|\le r_\theta$; thus choosing $\eta$ sufficiently small ensures, uniformly in $\theta$ and $\zeta$,
\[
 \sup_{\substack{|t|\le r_\theta\\|\zeta-\theta|\le1}}
 |H(t,\zeta)|\le C_0.
\]
By Cauchy's estimates in $t$, the Taylor coefficients
$H(t,\zeta)=\sum_{j\ge0}t^jP_j(\zeta)$ satisfy
$|P_j(\zeta)|\le C_0r_\theta^{-j}$ on this disk. More importantly,
for $|t|\le r_\theta/2$ the geometric-series tail gives
\begin{equation}\label{eq:schwartzcentralCauchy}
 \sup_{|\zeta-\theta|\le1}
 \left|H(t,\zeta)-\sum_{j=0}^2t^jP_j(\zeta)\right|
 \le 2C_0|t|^3r_\theta^{-3}
 \le C|t|^3(1+|\theta|)^6.
\end{equation}
The coefficients $P_1,P_2$ are precisely those of
Lemma~\ref{lem:global}: this follows by expanding $L(t\zeta)$ at $t=0$.
The analytic remainder to which the $\zeta$-Cauchy estimate is
applied is
\[
 \mathcal R_N(\zeta)=e^{ia\zeta}\bigl[H(N^{-1},\zeta)
 -1-N^{-1}P_1(\zeta)-N^{-2}P_2(\zeta)\bigr].
\]
Cauchy's derivative estimate on the concentric $\zeta$-disk of radius $1/2$, using the uniform bound on the disk of radius $1$,
then shows, for every integer $k\ge0$,
\begin{equation}\label{eq:schwartzcentralDerivatives}
 \max_{0\le j\le k}\left|\partial_\theta^j R_N(\theta)\right|
 \le C_k N^{-3}(1+|\theta|)^6,
 \quad |\theta|\le N^\varepsilon,
\end{equation}
where $R_N(\theta)=F_N(\theta,a)-e^{ia\theta}
\sum_{j=0}^2N^{-j}P_j(\theta)$, with the fixed choice $0<\varepsilon<1/2$ made above. Indeed, choose $N_0=N_0(\varepsilon,\eta,\delta)$ sufficiently
large that $N^{-1}\le r_\theta/2$ for all
$|\theta|\le N^\varepsilon$ and $N^\varepsilon<\delta N$
whenever $N\ge N_0$; this is possible because
$2\varepsilon<1$. Multiplication by
$e^{ia\zeta}$ introduces only a bounded factor on the complex disk.
The estimates hold with constants depending on $a,k,\varepsilon$,
but not on $N$ or $\theta$.

We next give the derivative bound needed outside this central region.
Choose $0<\delta<r_0/2$. For real $|\theta|\le\delta N$,
\[
 \partial_\theta^j[NL(\theta/N)]
 =N^{1-j}L^{(j)}(\theta/N),\qquad j\ge1.
\]
Every $L^{(j)}$ is bounded on $[-\delta,\delta]$. Consequently,
for each fixed $k\ge0$, the product rule (or the complete Bell
polynomial formula for derivatives of an exponential) yields
\begin{equation}\label{eq:schwartzintermediateDerivatives}
 |\partial_\theta^kF_N(\theta,a)|
 \le C_{k,\delta}|F_N(\theta,a)|
 \le C_{k,\delta}\exp\left(\frac{c\theta^2}{2N}\right),
 \quad |\theta|\le\delta N.
\end{equation}
Here the last inequality follows from the exact identity
$|F_N(\theta,a)|^2=(1+c\sin^2(\theta/N))^N$.
In particular the estimate is uniform in $N$; no factor $a^N$
occurs in this intermediate region. The derivatives of
$e^{ia\theta}\sum_{j=0}^2N^{-j}P_j(\theta)$ are bounded by
$C_k(1+|\theta|)^4$ for $N\ge1$.

Let $g(x)=\exp[-(x-x_0)^TQ(x-x_0)/2]$ as in
Lemma~\ref{lem:global}. There are constants $C_*,q_*>0$ such that
$|\partial_x^\gamma g(x)|\le C_\gamma(1+|x|)^{|\gamma|}
 e^{-q_*|x|^2}$ for every fixed $\gamma$.
On $N^\varepsilon<|\theta|\le\delta N$, the chain rule,
\eqref{eq:schwartzintermediateDerivatives}, and
$\theta^2\le |v|^2|x|^2$ show, for sufficiently large $N$, that
\[
 |x^\alpha\partial_x^\beta(g(x)R_N(v\cdot x))|
 \le C_{\alpha,\beta}(1+|x|)^{K_{\alpha,\beta}}
 e^{-q_*|x|^2/2}.
\]
Explicitly, fix $N_0$ such that
$\frac{(a^2-1)|v|^2}{2N}\le\frac{q_*}{2}$ for $N\ge N_0$;
this makes the preceding bound independent of $N$.
Because $|x|\ge N^\varepsilon/|v|$, Gaussian tail integration gives
an $L^2$ bound $O(e^{-c_1N^{2\varepsilon}})$ on this region.

On $|\theta|>\delta N$, differentiate the finite Fourier sum.
As $|s_{j,N}|\le1$ and $\sum_j|c_{j,N}(a)|=a^N$,
\[
 |\partial_\theta^kF_N(\theta,a)|\le a^N
 \quad(k\ge0,\ \theta\in\mathbb R).
\]
Since $|x|\ge\delta N/|v|$, multiplication by derivatives of $g$
and fixed polynomial weights gives an $L^2$ bound
$C_{\alpha,\beta}a^Ne^{-c_2N^2}=O(e^{-c_3N^2})$;
the polynomial correction terms obey the same bound without $a^N$.

Finally, on $|\theta|\le N^\varepsilon$ the estimates
\eqref{eq:schwartzcentralDerivatives}, the chain rule, and
Gaussian integrability imply
\[
 \|\mathbf1_{\{|\theta|\le N^\varepsilon\}}
 x^\alpha\partial_x^\beta(gR_N(v\cdot x))\|_2
 \le C_{\alpha,\beta}N^{-3}.
\]
Adding the three regions proves the result. The argument estimates
whole-space derivatives before restriction; it does not differentiate
the characteristic functions of the regions.
\end{proof}
\begin{theorem}[Weighted electric-field STFT remainder]\label{thm:weightedelectric}
Let $P:\R^d\to\R$ be a real polynomial, $b\in\R^d$, and
$U_{P,b}(t)=e^{-it[P(-i\nabla)-b\cdot x]}$.
Let $f_N$, $f_a$ be the fixed Gaussian directional signals of
Lemma~\ref{lem:global}, and let $P_1,P_2$ be its first two correction
polynomials. Put
\[
 r_j(x)=e^{ia v\cdot x}P_j(v\cdot x)g(x),\qquad j=1,2.
\]
For any $m\in\N_0$, $T>0$, and $h\in\mathcal S(\R^d)$,
\[
 \sup_{|t|\le T}
 \left\|V_hU_{P,b}(t)
 \left(f_N-f_a-N^{-1}r_1-N^{-2}r_2\right)
 \right\|_{L^2((1+|u|+|\omega|)^{2m}du\,d\omega)}
 =O_{m,T}(N^{-3}).
\]
The implicit constant may depend on $m,T,P,b,h,a,v,Q,x_0$,
and on the fixed choice of norms, but not on $N$. The number of
Schwartz seminorms required depends on $\deg P$.
The motivating separable even-polynomial momentum Hamiltonians are
included as a special case.
\end{theorem}
\begin{proof}
Let $E_N$ be the initial-data remainder. Lemma~\ref{lem:schwartzremainder} gives
$\|x^\alpha\partial^\beta E_N\|_2\le C_{\alpha,\beta}N^{-3}$
for every pair of multi-indices. By the Fourier-characteristic formula~\eqref{fourierprop},
$\widehat{U_{P,b}(t)E_N}$ is a translate of $\widehat E_N$ multiplied by
$e^{-i\Theta_t(\xi-tb)}$. On $|t|\le T$, every fixed derivative of this phase is
bounded by a polynomial in $\xi$, of degree depending only on $P$ and the differentiation order, with coefficients uniform in $t$. Leibniz's rule, Plancherel, and the Fourier identities relating $x^\alpha\partial^\beta$ to frequency derivatives and powers therefore imply, for each $M$,
\[
\sup_{|t|\le T}\sum_{|\alpha|+|\beta|\le M}
\|x^\alpha\partial^\beta U_{P,b}(t)E_N\|_2
\le C_{M,T}\sum_{|\gamma|+|\delta|\le M'}
\|x^\gamma\partial^\delta E_N\|_2=O_{M,T}(N^{-3}),
\]
where one may take $M'=(r+1)(M+1)$ with
$r=\max\{1,\deg P\}$. To justify this finite loss explicitly,
write the Fourier-side propagator as
\[
 \widehat{U_{P,b}(t)E_N}(\xi)
 = e^{-i\Theta_t(\xi-tb)}\widehat E_N(\xi-tb).
\]
For $|t|\le T$ and each multi-index $\mu$, the chain rule
(Fa\`a di Bruno formula) expresses
$\partial_\xi^\mu e^{-i\Theta_t(\xi-tb)}$ as this
unit-modulus exponential times a polynomial in $\xi$ of degree
at most $(r-1)|\mu|$, whose coefficients are bounded uniformly for
$|t|\le T$. Consequently, differentiating $\alpha$ times and
multiplying by $\xi^\beta$ produces finitely many terms of the form
\[
 p_{\alpha,\beta,\gamma,t}(\xi)
 \bigl(\partial_\xi^\gamma\widehat E_N\bigr)(\xi-tb),
 \qquad |\gamma|\le|\alpha|,
\]
where $\deg p_{\alpha,\beta,\gamma,t}\le
|\beta|+(r-1)|\alpha|$. Translation $\zeta=\xi-tb$ changes
these polynomials into polynomials in $\zeta$ with uniformly bounded
coefficients on $|t|\le T$. Applying Plancherel, recall explicitly that multiplication by
$\zeta^\nu$ on the Fourier side corresponds to $D^\nu$ in physical
space, while $\partial_\zeta^\gamma$ corresponds (up to powers of
$i$) to multiplication by $x^\gamma$. Hence Plancherel bounds each
term by finitely many norms
$\|x^\gamma D^\delta E_N\|_2$ with
$|\gamma|+|\delta|\le rM$ when
$|\alpha|+|\beta|\le M$. The larger displayed choice of $M'$
therefore suffices. All constants depend on the fixed parameters
$M,T,P,b$ and the Fourier convention, not on $N$.
Apply Lemma~\ref{lem:weighted} to conclude.
\end{proof}
\begin{remark}
This result is a polynomially weighted STFT statement on every bounded
time interval. It does not claim concentration in a moving tube around
the nonlinear Hamiltonian trajectory or bounds uniform as $T\to\infty$.
Those refinements require additional microlocal estimates.
\end{remark}

\section{Exact co-moving observables for electric-field polynomial Hamiltonians}\label{sec:comoving}
The preceding weighted STFT estimates can be expressed in an exact, time-dependent operator norm. The resulting identities are useful for analyzing polynomial Hamiltonians of any degree, but should not be confused with an STFT localization estimate near a scalar classical trajectory.
\begin{theorem}[Exact conjugation of canonical observables]\label{thm:comoving}
Let $P:\R^d\to\R$ be a real polynomial, $b\in\R^d$, $D=-i\nabla$, and $U(t)=e^{-it(P(D)-b\cdot x)}$. Define
\[
 \Theta_t(\xi)=\int_0^t P(\xi+sb)\,ds,
 \qquad
 X_t=x-\nabla\Theta_t(D-tb),\qquad \Xi_t=D-tb.
\]
On $\mathcal S(\R^d)$ one has the exact identities
\[
 X_tU(t)=U(t)x,\qquad \Xi_tU(t)=U(t)D.
\]
Consequently, for any multi-indices $\alpha,\beta$,
\[
 \|X_t^\alpha\Xi_t^\beta U(t)f\|_2
 =\|x^\alpha D^\beta f\|_2.
\]
Thus the co-moving graph norms
\[
 \|u\|_{m,t}:=\sum_{|\alpha|+|\beta|\le m}\|X_t^\alpha\Xi_t^\beta u\|_2
\]
are exactly preserved by evolution: $\|U(t)f\|_{m,t}=\|f\|_{m,0}$. In particular, an expansion through order $N^{-2}$ with a
Schwartz graph-norm remainder $O(N^{-3})$ is propagated with the
same $O(N^{-3})$ remainder in co-moving graph norms, uniformly for
all $t\in\R$.
\end{theorem}
\begin{proof}
The Fourier characteristic formula gives
\[
 U(t)=M_{tb}\,e^{-i\Theta_t(D)}.
\]
The modulation identity $M_{tb}^{-1}(D-tb)M_{tb}=D$ implies
$U(t)^{-1}\Xi_tU(t)=D$. For the position operator, use the Fourier representation $x_j=i\partial_{\xi_j}$ to compute
\[
 e^{i\Theta_t(D)}x_j e^{-i\Theta_t(D)}
 =x_j+\partial_j\Theta_t(D).
\]
It follows that $U(t)^{-1}x_jU(t)=x_j+\partial_j\Theta_t(D)$, while
$U(t)^{-1}\partial_j\Theta_t(D-tb)U(t)=\partial_j\Theta_t(D)$.
Subtracting yields $U(t)^{-1}X_{t,j}U(t)=x_j$. The identities for ordered monomials follow by composition on the invariant Schwartz core. Finally, unitarity of $U(t)$ gives the norm equality. Schwartz invariance follows directly from the Fourier characteristic formula, because all derivatives of the polynomial phase grow at most polynomially and modulation and translation preserve Schwartz functions. The identities hold first on the common invariant Schwartz core.
For a fixed ordered monomial $A_{\alpha,\beta}=x^\alpha D^\beta$
with its closure from $\mathcal S(\R^d)$, define the co-moving
monomial as the closed operator
\[
 A_{\alpha,\beta}(t):=U(t)A_{\alpha,\beta}U(t)^{-1},
 \qquad
 \operatorname{Dom}A_{\alpha,\beta}(t)
 =U(t)\operatorname{Dom}A_{\alpha,\beta}.
\]
On $\mathcal S$ this agrees with the ordered product
$X_t^\alpha\Xi_t^\beta$. Thus the closed monomial identity and
norm equality follow from unitary conjugation, without assuming
automatic domain identities for arbitrary products of unbounded
operators.
\end{proof}
\begin{remark}
For nonlinear $P$, $X_t$ is an operator-valued shift involving the momentum observable. This theorem gives an exact quantum co-moving framework for weighted estimates; an estimate localizing the ordinary STFT near a scalar Hamiltonian trajectory requires further semiclassical or narrow-band hypotheses.
\end{remark}

\begin{remark}[Scope of the weighted propagation result]
The constants in Theorem~\ref{thm:weightedelectric} depend on the fixed
Gaussian parameters, the fixed STFT window, the time horizon $T$, the
weight order $m$, and the coefficients of $P$ and $b$. They are independent
of $N$ and $t\in[-T,T]$. The estimate is a finite-time polynomially weighted
bound, not a uniform-in-time assertion in ordinary weighted STFT spaces.
\end{remark}

\section{Coupled polynomial frequency maps: the remaining analytic step}
For $\Lambda(s)=(s^{p_1},\dots,s^{p_d})$ the expression
\[
 \sum_{j=0}^N c_{j,N}(a)e^{i\Lambda(s_{j,N})\cdot x}
\]
is not a product of independent one-dimensional supershifts. A formal
Taylor expansion in $s$ suggests the correction
\[
 \frac{1-a^2}{2N}
 \big(i\Lambda''(a)\cdot x-(\Lambda'(a)\cdot x)^2\big)
 e^{i\Lambda(a)\cdot x},
\]
but the signed coefficients have total variation $a^N$.
Consequently a pointwise Taylor remainder bounded uniformly in $j$ is
insufficient to sum the weighted errors. The multivariable analytic
supershift literature \cite{CPSS,CSSY} supplies qualitative convergence
under suitable holomorphic assumptions. A fully quantitative
$O(N^{-2})$ remainder for the coupled curve, uniform in the weighted
window norms required for frame asymptotics, is retained as a separate
problem, rather than asserted without proof. The directional case in
Sections 3--5 is a rigorous benchmark for this extension.

\section{Conclusions and limitations}
The results provide rigorous quantitative supershift asymptotics in global STFT norms, compact-window expansions for directional frequency families, perturbative Gabor reconstruction, precise Gaussian entire-growth membership, and electric-field evolution. The positive fixed-step Gaussian theorem and the critical-type divergence result must be read together: membership of individual kernels is weaker than convergence, and the endpoint can fail dramatically. The widening-Gaussian theorem gives a sharp global $L^2$ energy transition and a critical rescaled intensity profile, but not a moving-window STFT profile. An independent scalar Gaussian reduction verifies the critical leading constant and identifies the rank-one nature of the degeneracy. Nonlinear Hamiltonian-flow localization and local critical-window STFT profiles remain further problems. The elementary Hilbert-space persistence results should not be confused with such genuinely new estimates.

\end{document}